\documentclass[11pt,a4paper]{amsart}
\usepackage{amsmath,amsthm,amssymb}
\usepackage{graphicx}
\usepackage{comment}
\newcommand{\mychoice}[3]{#1
}
\mychoice{
\newcommand{\plabel}[1]{ \label{#1}}
\newcommand{\gbibitem}[1]{ \bibitem{#1}}
\newcommand{\snewpage}{}

\specialcomment{commentx}{}{}\excludecomment{commentx}
\specialcomment{commenty}{}{}\excludecomment{commenty}
\specialcomment{commentz}{}{}

}{
\usepackage{xcolor}
\newcommand{\plabel}[1]{ \label{#1}\rlap{\smash{${}^{^{[#1]}}$}}}
\newcommand{\gbibitem}[1]{ \bibitem{#1}\rlap{\smash{${}^{^{[#1]}}$}}}
\newcommand{\snewpage}{\newpage}

\newenvironment{commentx}{\color{magenta} }{\color{black} }
\newenvironment{commenty}{\color{blue} }{\color{black} }

}{
\usepackage{xcolor}
\newcommand{\plabel}[1]{ \label{#1}}
\newcommand{\gbibitem}[1]{ \bibitem{#1}}
\newcommand{\snewpage}{}

\specialcomment{commentx}{}{}\excludecomment{commentx}

\specialcomment{commentz}{}{}\excludecomment{commentz}

}
\DeclareMathOperator{\sgn}{sgn}

\DeclareMathOperator{\Lie}{Lie}
\DeclareMathOperator{\erf}{erf}
\DeclareMathOperator{\proj}{pr}\DeclareMathOperator{\gen}{gen}
\DeclareMathOperator{\gr}{\boldsymbol{\mathbf{gr}}}

\DeclareMathOperator{\ad}{ad}

\DeclareMathOperator{\dist}{dist}

\newcommand{\leaveout}[1]{}

\theoremstyle{definition}
\newtheorem{point}{}[section]
\newtheorem{remark}[point]{Remark}

\theoremstyle{plain}
\newtheorem{prop}[point]{Proposition}

\newtheorem{theorem}[point]{Theorem}

\newcommand{\qedremark}{  \renewcommand{\qedsymbol}{$\triangle$} \qed \renewcommand{\qedsymbol}{$\Box$}}

\newcommand{\marginextend}[1]{ \addtolength{\oddsidemargin}{-#1}  \addtolength{\evensidemargin}{-#1}\addtolength{\textwidth}{#1}\addtolength{\textwidth}{#1}}
\newcommand{\updownextend}[1]{ \addtolength{\topmargin}{-#1}  \addtolength{\textheight}{#1}
\addtolength{\textheight}{#1}}
\marginextend{1.25cm}
\updownextend{0cm}
\allowdisplaybreaks[3]
\begin{document}
\title{The noncommutative heat equation and certain Lie series}
\author{Gyula Lakos}
\email{lakos@renyi.hu}
\address{Alfréd Rényi Institute of Mathematics}
\begin{abstract}
We approach the convergence of the Magnus, Wilcox, and symmetric Wilcox expansions
 by a non-commutative heat equation derived from the Maurer--Cartan equation.
\end{abstract}
\maketitle

\section{Introduction}\plabel{sec:intro}
\textbf{Lie series, Maurer--Cartan equation, heat equation.}
The Magnus expansion is a continuous generalization of the ``discrete'' Baker--Campbell--Hausdorff expansion.
It is due to Magnus \cite{Mag}, although it has been rediscovered a few times.
A classical review of the Magnus expansion is given by Blanes, Casas, Oteo, Ros \cite{BCOR},
 which also contains some information regarding other Lie expansions.
More recent information in those directions is given by Arnal,  Casas,  Chiralt, Oteo  \cite{ACCO}.

The objective of this paper is to provide intuitive although not particularly sharp convergence bounds for the
Magnus expansion, and for the related Wilcox expansion, and also for versions of the symmetric Wilcox expansion.
Our approach uses the Maurer--Cartan equation  and the associated non-commutative heat equation.
Intuitive is, however, not the same as technically unassuming.
Thus, for the sake of ease of presentation, we will pursue a rather relaxed style in the discussion, in which
 several technical estimates are omitted, but which can be made completely precise by the professional mathematician.

\textbf{On the basic setting.}
The expansions above including the associated exponential formulae are the easiest to be understood in Banach algebraic sense.
However, as we deal with Lie polynomials, the norm estimates for the terms of the expansions can also be done Banach--Lie algebraic sense.
But then the meanings of the exponential formulae are not clear.
They either make sense in the adjoint representation or in an appropriate setting of Banach--Lie groups.
Nevertheless, then these estimates can also be transferred to the setting of Lie groups.
For this reason, and also because the Banach--Lie algebraic estimates can also be used for the Banach algebraic setting
 (although not that effectively), we will primarily develop Banach--Lie algebraic estimates here.

\textbf{Outline of content.}
In Section  \ref{sec:lieser}, we review the Lie series we will deal with, and also some associated terminology.
In Section  \ref{sec:mcheat}, we try to  motivate and understand the forthcoming developments in the
 Lie group theoretic setting, which might be more familiar to many.
We consider the Maurer--Cartan equation; and, more specifically, we choose the heat (diffusion) prescription to drive it.
In Section \ref{sec:formheat}, we return to the Banach and Banach--Lie algebraic viewpoints.
We set up (formal) solutions to the non-commutative heat equation directly, and find out what ``trivial'' estimates follow from it.
In Section \ref{sec:perex}, we work out an example of measures over $2\times 2$ real matrices in the periodic setting.
\snewpage

\section{Preliminaries: Some Lie expansions}\plabel{sec:lieser}
\textbf{The Magnus expansion.}
Let us recall some information on the Magnus expansion (series) in the spirit of \cite{L1}.
The simplest, formal setting of the Magnus expansion is as follows.
Assume that $\phi$ is a Banach algebra $\mathfrak A$ valued continuous measure of finite variation on a possibly infinite interval $I$.
Let $T$ be a formal ``commutative'' variable.
Then the time-ordered exponential of $T\cdot\phi$ is defined as
\[\exp_{\mathrm R}(T\cdot\phi)=1+\sum_{n=1}^\infty T^n\cdot\int_{t_1\leq\ldots\leq t_n\in I} \phi(t_1) \ldots \phi(t_n). \]
In this case,
\begin{equation}
\exp_{\mathrm R}(T\cdot\phi)=\exp(\mu_{\mathrm R}(T\cdot\phi))
\plabel{eq:expmag}
\end{equation}
with
\begin{equation}
\mu_{\mathrm R}(T\cdot\phi)=\sum_{n=1}^\infty T^n\mu_{\mathrm R,n}(\phi)
\plabel{eq:magexp}
\end{equation}
such that
\[\mu_{\mathrm R,n}(\phi)=\int_{t_1\leq\ldots\leq t_n\in I}\mu_n(\phi(t_1),\ldots,\phi(t_n)),\]
where $\mu_n(X_1,\ldots,X_n)$ are commutator polynomials which are linear in their variables.
For example,
\begin{align*}
\mu_1(X_1)&=X_1,\\
\mu_2(X_1,X_2)&=\frac12[X_1,X_2],\\
\mu_3(X_1,X_2,X_3)&=\frac16[[X_1,X_2],X_3]+\frac16[X_1,[X_2,X_3]],\\
\mu_4(X_1,X_2,X_3,X_4)&=\frac1{12}[[X_1,[X_2,X_3],X_4]+\frac1{12}[X_1,[[X_2,X_3],X_4]]+\\
&\quad+\frac1{12}[[X_1,X_2],[X_3,X_4]]+\frac1{12}[[X_1,X_3],[X_2,X_4]].
\end{align*}
Here \eqref{eq:magexp} is the Magnus series, and \eqref{eq:expmag} is the associated exponential identity which holds for algebraic reasons.
The Magnus expansion in the classical (original) setting is when we substitute $T=1$.
This is, in general, problematic, but it can be done if
 \begin{equation}
\sum_{n=1}^\infty |\mu_{\mathrm R,n}(\phi)|_{\mathfrak A}<+\infty.
\plabel{eq:magabs}
\end{equation}
This is the case of (absolute) convergence of the Magnus series.
Using combinatorial arguments, it can be shown that
\begin{equation}
\sum_{n=1}^\infty|\mu_{\mathrm R,n}(\phi)|_{\mathfrak A}T^n
\stackrel{\forall T^n}{\leq}\sum_{n=1}^\infty\frac1{2^{n-1}} \left( \textstyle{\int|\phi|_{\mathfrak A}}\right)^n  T^n.
\plabel{eq:magmo}
\end{equation}
(The notation means that the relation holds in the coefficients of the powers of $T$.)

Consequently,
if $\int|\phi|_{\mathfrak A}<2$ (a bound for the associated cumulative norm) holds, then
\[\sum_{n=1}^\infty|\mu_{\mathrm R,n}(\phi)|_{\mathfrak A}\leq
\frac{\left( \textstyle{\int|\phi|_{\mathfrak A}}\right)}{1-\frac12\left( \textstyle{\int|\phi|_{\mathfrak A}}\right)} <+\infty,\]
thus absolute convergence holds.
This is essentially the result of Moan Oteo \cite{MO}.
Although \eqref{eq:magmo} itself is not sharp, the associated cumulative norm radius $2$ is.

As the $\mu_n$ are commutator polynomials,
one can also obtain estimates for the terms of the Magnus expansion with respect to a
Lie algebra $\mathfrak g$ endowed with a  Banach--Lie norm $\|\cdot\|_{\mathfrak g}$
(where $\|[X,Y]\|_{\mathfrak g}\leq\|X\|_{\mathfrak g},\|Y\|_{\mathfrak g}$).
See \cite{L3} for the convergence norm estimate in the Banach--Lie case, where
the cumulative radius can be improved to more than $2.4$ .
In fact, using appropriate restrictions in the norm, and analytic continuation, we can obtain results in the setting of Lie groups
(see the next section).


\textbf{The Wilcox expansion.}
A similar expansion is the ``Magnus--Zassenhaus'' or Wilcox expansion.
(In what follows: `Wilcox expansion'.)
In the formal case, it is such that
\[\exp_{\mathrm R}(T\cdot \phi)=
\ldots
\exp(T^3\cdot\zeta^\leftarrow_{\mathrm R,3}(\phi))
\exp(T^2\cdot\zeta^\leftarrow_{\mathrm R,2}(\phi))
\exp(T\cdot\zeta^\leftarrow_{\mathrm R,1}(\phi))
\]
with
\[\zeta^\leftarrow_{\mathrm R,n}(\phi)=\int_{t_1\leq\ldots\leq t_n\in I}\zeta_n^+(\phi(t_1),\ldots,\phi(t_n)),\]
where $\zeta^\leftarrow_n(X_1,\ldots,X_n)$ are again  commutator  polynomials which are linear in their variables.
For example,
\begin{align*}
\zeta^\leftarrow_1(X_1)&=X_1,\\
\zeta^\leftarrow_2(X_1,X_2)&=\frac12[X_1,X_2],\\
\zeta^\leftarrow_3(X_1,X_2,X_3)&=\frac13[X_2,[X_1,X_3]]+\frac13[X_1,[X_2,X_3]], \\
\zeta^\leftarrow_4(X_1,X_2,X_3,X_4)&=\frac14[X_1,[X_2,[X_3,X_4]]]+\frac14[X_1,[X_3,[X_2,X_4]]]+\frac14[X_2,[X_3,[X_1,X_4]]] .
\end{align*}
Now, the classical case is when we put $1$ to the place $T$.
Again, this works out if
 \begin{equation}
\sum_{n=1}^\infty |\zeta^\leftarrow_{\mathrm R,n}(\phi)|_{\mathfrak A}<+\infty.
\plabel{eq:wilabs}
\end{equation}
This expansion appears in BCH type form in Magnus \cite{Mag} as the ``Zassenhaus formula'',
 and in continuous form in Wilcox \cite{Wil} (calling it erroneously Fer's expansion, but which is, in fact, a different thing).

A variant of the expansion above is given as
\[\exp_{\mathrm R}(T \cdot\phi)=
\exp(T \cdot\zeta^\rightarrow_{\mathrm R,1}(\phi))
\exp(T^2\cdot\zeta^\rightarrow_{\mathrm R,2}(\phi))
\exp(T^3\cdot\zeta^\rightarrow_{\mathrm R,3}(\phi))
\ldots
\]
where
\[\zeta^\rightarrow_{\mathrm R,n}(\phi)=\int_{t_1\leq\ldots\leq t_n\in I}\zeta_n^\rightarrow(\phi(t_1),\ldots,\phi(t_n)),\]
and
\[\zeta^\rightarrow_n(X_1,\ldots,X_n)=-\zeta^\leftarrow_n(-X_n,\ldots,-X_1).\]
This version $\zeta^\rightarrow$ is completely analogous and, in fact, equivalent to $\zeta^\leftarrow$
(by passing to the negative transposed measure).

Convergence has mainly been investigated in the Banach algebraic setting.
In the case of the Zassenhaus formula, after some simpler estimates by Suzuki \cite{Suz},
Bayen \cite{Bay} proves convergence when the cumulative norm $\int|\phi|_{\mathfrak A}$ is less than $0.5967\ldots$ (a solution of an equation).
Numerical results by Casas, Murua, Nadinic \cite{CMN}, however, suggest convergence for cumulative norm less than $1.054$ .
In the case of the Wilcox expansion, numerical results by  Arnal,  Casas,  Chiralt, Oteo  \cite{ACCO}
suggest convergence for cumulative norm less than $0.6584$ .
(The Banach algebraic bounds also apply in the Banach--Lie case.)
\snewpage

\textbf{Symmetric Wilcox expansions.}
Another variant is the inward expanding symmetric Wilcox expansion.
In the formal case, it is such that
\begin{multline*}
\exp_{\mathrm R}(T\cdot \phi)=\exp(T\cdot\eta^{\bowtie}_{\mathrm R,1}(\phi)/2)
\exp(T^2\cdot\eta^{\bowtie}_{\mathrm R,2}(\phi)/2)\exp(T^3\cdot\eta^{\bowtie}_{\mathrm R,3}(\phi)/2)\ldots\cdot
\\
\cdot\ldots\exp(T^3\cdot\eta^{\bowtie}_{\mathrm R,3}(\phi)/2)\exp(T^2\cdot\eta^{\bowtie}_{\mathrm R,2}(\phi)/2)\exp(T\cdot\eta^{\bowtie}_{\mathrm R,1}(\phi)/2)
\end{multline*}
with
\[\eta^{\bowtie}_{\mathrm R,n}(\phi)=\int_{t_1\leq\ldots\leq t_n\in I}\eta^{\bowtie}_n(\phi(t_1),\ldots,\phi(t_n)),\]
where $\eta^{\bowtie}_n(X_1,\ldots,X_n)$ are again (commutator) polynomials which are linear in their variables.
For example,
\begin{align*}
\eta^{\bowtie}_1(X_1)&=X_1,\\
\eta^{\bowtie}_2(X_1,X_2)&=\frac12[X_1,X_2],\\
\eta^{\bowtie}_3(X_1,X_2,X_3)&= \frac16[[X_1,X_2],X_3]+\frac16[X_1,[X_2,X_3]],\\
\eta^{\bowtie}_4(X_1,X_2,X_3,X_4)&= \frac18[[[X_1,X_2],X_3],X_4]+\frac18[X_1,[X_2,[X_3,X_4]]].
\end{align*}
Similarly as before, we can put $1$ to the place $T$ if
 \begin{equation}
\sum_{n=1}^\infty |\eta^{\bowtie}_{\mathrm R,n}(\phi)|_{\mathfrak A}<+\infty.
\plabel{eq:swilabs}
\end{equation}

Yet another variant is the outward expanding symmetric Wilcox expansion.
In the formal case, it is such that
\begin{multline*}
\exp_{\mathrm R}(T^3\cdot \phi)=\ldots\exp(T\cdot\eta^{\leftrightarrow}_{\mathrm R,3}(\phi)/2)
\exp(T^2\cdot\eta^{\leftrightarrow}_{\mathrm R,2}(\phi)/2)\exp(T\cdot\eta^{\leftrightarrow}_{\mathrm R,1}(\phi)/2)\cdot
\\
\cdot\exp(T\cdot\eta^{\leftrightarrow}_{\mathrm R,1}(\phi)/2)\exp(T^2\cdot\eta^{\leftrightarrow}_{\mathrm R,2}(\phi)/2)\exp(T^3\cdot\eta^{\leftrightarrow}_{\mathrm R,3}(\phi)/2)\ldots
\end{multline*}
with
\[\eta^{\leftrightarrow}_{\mathrm R,n}(\phi)=\int_{t_1\leq\ldots\leq t_n\in I}\eta^{\leftrightarrow}_n(\phi(t_1),\ldots,\phi(t_n)),\]
where $\eta^{\leftrightarrow}_n(X_1,\ldots,X_n)$ are again (commutator) polynomials which are linear in their variables.
For example,
\begin{align*}
\eta^{\leftrightarrow}_1(X_1)&=X_1,\\
\eta^{\leftrightarrow}_2(X_1,X_2)&=\frac12[X_1,X_2],\\
\eta^{\leftrightarrow}_3(X_1,X_2,X_3)&= \frac16[[X_1,X_2],X_3]+\frac16[X_1,[X_2,X_3]],\\
\eta^{\leftrightarrow}_4(X_1,X_2,X_3,X_4)&= \frac18\left([X_1,[X_4,[X_3,X_2]]] +\frac18[X_2,[X_3,[X_4,X_1]]]\right)\\
&\quad+ \frac18\left( [X_3,[X_2,[X_4,X_1]]]+\frac18[X_4,[X_1,[X_3,X_2]]]\right)
\end{align*}
(with the two last summands actually being equal).
Similar remark applies to convergence.
The inward and outward expanding symmetrical Wilcox expansions are not completely analogous.

Applied to the BCH measure $X\mathbf 1_{[0,1)}\boldsymbol.Y\mathbf 1_{[1,2)}$, we obtain
 the BCH type forms of these expansions, which are the inward and outward expanding ``symmetric Zassenhaus formulas'', respectively.
A kind of symmetric Zassenhaus formula but which is not the expansion of the
 time-ordered exponential is examined in Arnal, Casas, Chiralt \cite{ACC}.
Otherwise, the convergence of these symmetric expansions seems not have been investigated particularly.

\textbf{The unicity of the expansions in the formal case.}
The grading according to $T$  allows the reconstruction of the Lie series mentioned above from the time-ordered exponentials.
For example, in case of $\eta^\leftrightarrow$, assume that
\begin{multline*}
\exp_{\mathrm R}(T\cdot\alpha)=\ldots (\exp T^3\cdot A_3/2)(\exp  T^2\cdot A_2/2)(\exp T^1\cdot A_1/2)
\\\cdot
(\exp  T^1\cdot A_1/2)(\exp T^2\cdot A_2/2)(\exp T^3\cdot A_3/2)\ldots,
\end{multline*}
 where $A_i$ are elements from our Banach algebra.
Then the $A_i$ can be obtained by the recursion
\begin{multline*}
A_i = \text{the coefficient of $T^i$ in}
\\ \log\biggl( \left((\exp T^{i-1}\cdot A_{i-1}/2)\ldots(\exp  T^2\cdot A_2/2)(\exp T^1\cdot A_1/2))\right)^{-1} \\
\cdot\exp_{\mathrm R}(\mathrm Z^1_{\nu})\cdot\left( (\exp  T^1\cdot A_1/2)(\exp T^2\cdot A_2/2)\ldots(\exp T^{i-1}\cdot A_{i-1}/2) \right)^{-1} \biggr).
\end{multline*}
The expression above can be resolved by multiple Magnus and BCH expansions, thus we know that
 the $A_i$ can be obtained as integrals of commutators.
The actual shape of the Lie polynomials involved can be tested by multiple BCH expansions.
Similar comment applies to the other Lie series.
In general, the entire function $\lambda\in\mathbb C\mapsto \exp_{\mathrm R}(\lambda\cdot\alpha)$,
 the formal time-ordered exponential $\exp_{\mathrm R}(T\cdot\alpha)$, the Magnus series, the left and
 right expanding Wilcox series, the inward  and outward expanding symmetric Wilcox series of $\alpha$ contain
 the same information; and one can be computed from the other.
A useful viewpoint is provided by
\snewpage

\textbf{Formal noncommutative masses.}
In what follows all measures will be induced by intervals, the base algebra can be $\mathbb R$ or $\mathbb C$.
In \cite{L1} and \cite{L3} we have introduced `tautological' measures $\mathrm Z^1_{[a,b)}$ (Banach algebraic) and
 $\mathrm Z^{1,\Lie}_{[a,b)}$ (Banach--Lie algebraic).
For example $\mathrm Z^1_{[a,b)}$ is a Banach algebra $\mathrm F^{1}([a,b))$ valued measure
 such that for any subinterval  $J\subset [a,b)$, the inequality
 $|\mathrm Z^1_{[a,b)}(J)|_{\mathrm F^{1}([a,b))}\leq|J|$ (in fact: equality) holds;
 and the Banach algebra $\mathrm F^{1}([a,b))$ is generated by  the
 $ \mathrm Z^1_{[a,b)}(J)$, with the largest possible norms.
Similar construction applies in the Banach--Lie case.
These are ``free non-commutative mass'' versions of the Lebesgue measure.
Analogues for other variation measures can be defined similarly, but they can also be defined as reparametrizations.
For example, in the Banach algebraic case:
Assume that $\nu$ is a nonnegative measure on the interval $I$.
Then we may define $\mathrm Z^1_{\nu}$ by
\[\mathrm Z^1_{\nu}(\{t\})=\mathrm Z^1_{[0,\int\nu)}( [\nu(\{x\in I:x<t\}),\nu(\{x\in I:x\leq t\}))  ),\]
\[\mathrm Z^1_{\nu}(\{(t_1,t_2)\})=\mathrm Z^1_{[0,\int\nu)}( [\nu(\{x\in I:x\leq t_1\}),\nu(\{x\in I:x< t_2\}))  ).\]
Thus  $\mathrm Z^1_{\nu}$ is just a reparametrized version of $\mathrm Z^1_{[0,\int\nu)}$ but with variation measure $\nu$.
(The actually generated ambient algebra is smaller than  $\mathrm F^{1}([a,b))$  for a non-continuous measure.)

If $\alpha$ is a Banach algebra $\mathfrak A$ valued measure with finite variation, then
 there is a natural contractive homomorphism from $ \mathrm F^1([0,\int\nu)) $
 (the sort of free Banach algebra generated by the values of the measure) into   $\mathfrak A$,
 taking $\mathrm Z^1_{|\alpha|_{\mathfrak A}}$ into $\alpha$.
Here the role of $\mathrm Z^1_{|\alpha|_{\mathfrak A}}$ is analogous to a collection of variables for polynomials,
 both universal and combinatorially recognizable.
In general, if norm estimates allow it, in computations, $\alpha$ can be replaced by  $\mathrm Z^1_{|\alpha|_{\mathfrak A}}$,
 and then the contractive homomorphism can be taken.
Note that $ \mathrm F^1([0,\int\nu)) $ is a naturally graded Lie algebra where the values of the measure
 $\mathrm Z^1_{\nu}$ are of uniformly grade $1$
 (and this also allows to use the ambient locally convex algebra $\mathrm F^{1,\mathrm{loc}}([0,\int\nu))$).
In other terms, there is a natural way to scaling by $T$.
The formal reconstruction regarding our Lie series explained earlier can be applied to $\mathrm Z^1_{\nu}$.
This recursion can actually be used to define the expansion in universal ``non-commutative polynomial terms''.

A  similar discussion applies in the Banach--Lie algebraic case.
There the exponential formulae do not (necessarily) make sense but the terms of the series
 can be evaluated and possibly estimated in the Banach--Lie setting.
The ambient Banach--Lie algebra $\mathrm F^{1,\Lie}([0,\int\nu))$ of $\mathrm Z^{1,\Lie}_{\nu}$
 is contracted relative to  the ambient Banach  algebra $\mathrm F^{1 }([0,\int\nu))$ of $\mathrm Z^{1 }_{\nu}$
 (restricting  $\mathrm F^{1 }([0,\int\nu))$ or extending $\mathrm F^{1,\Lie}([0,\int\nu))$ appropriately, cf. \cite{L3}).
However, on Lie expressions of a fixed degree they are comparable.
(If the degree is $k$, then at most by a factor $2^{k-1}$.)
Thus either of them can be used for formal solutions,
 but the formal Banach--Lie algebraic convergence estimates are typically better but by a factor at most $2$ in convergence radii.
(Doubling the Banach algebra norm yields a Banach--Lie norm;
while $ \mathrm F^{1,\Lie}([0,\int\nu))$ can be extended to a Banach algebra, contracted from $\mathrm F^{1 }([0,\int\nu))$.)

Our stances therefore will be somewhat peculiar:
When we are to obtain convergence estimates on Lie series, we always implicitly use formal noncommutative masses.
For the sake of arithmetic we use the Banach algebraic picture, but then we do the actual estimates in the Banach--Lie algebraic setting.

\textbf{On the non-Magnus norm estimates.}
Our objective is to establish baseline, trivial estimates which are trivial in the sense that very little combinatorics is used.
Nevertheless, our approach has the possibility to obtain much stronger estimates, although those would involve
 much more difficult computations.

\snewpage
\section{The Maurer--Cartan equation and a noncommutative heat equation}\plabel{sec:mcheat}
\textbf{The setting of Lie groups.}
In the viewpoint of Lie groups, the Magnus expansion can be imagined as follows.
If a, say, smooth path $g:[0,1]\rightarrow G$ into a Lie group is given such that
 $g(0)=1$ and $g(0)=h$, then, under favourable circumstances, the Magnus expansion  computes a value  $H$ such that $h=\exp H$.
It does this in the form of a sum of time-ordered integrals of higher commutators of
 $\alpha( x )=g( x )^{-1}\frac{\mathrm d}{\mathrm dt}g( x )$.
The expression for the Magnus expansion is
\[\sum_{n=1}^\infty\mu_{\mathrm R,n}(\alpha)\equiv
\sum_{n=1}^\infty\int_{0\leq  x _1\leq\ldots\leq  x _n\leq1}\mu_n(\alpha( x _1),\ldots,\alpha( x _n))
\,\mathrm d x _1\ldots \mathrm d x _n.\]
It can be said that the Magnus expansion straightens the Lie development $\alpha:[0,1]\rightarrow\mathfrak g$
 into a constant development given by $ x \in[0,1]\mapsto H$ (but which keeps the endpoints of the ``integrated'' path).
I.~e.~the sum $\mu_{\mathrm R}(\alpha)=\sum_{n=1}^\infty\mu_{\mathrm R,n}(\alpha)$
exists, and $\exp\mu_{\mathrm R}(\alpha)=g(0)^{-1}g(1)=h $.
As Lie groups are analytical objects, ``favourable circumstances'' include the
 case when the convergence radius of the power series
\[\lambda\in\mathbb C\mapsto \sum_{n=1}^\infty\lambda^n\mu_{\mathrm R,n}(\alpha)\]
 around $\lambda=0$ is greater than $1$.
(We have convergence and the exponential identity valid for $\lambda\cdot\alpha$ by Ado's theorem for $\lambda\sim0$.
Complexification is harmless here.
Then, having meaningful terms to compare in the exponential identity, the uniqueness of analytical continuation applies along $\lambda\in[0,1]$.)
Conditions for this can be provided if $\mathfrak g$ is endowed with a Banach--Lie norm $\|\cdot\|_{\mathfrak g}$
 (where $\|[X,Y]\|_{\mathfrak g}\leq\|X\|_{\mathfrak g}\,\|Y\|_{\mathfrak g}$).
For example, the convergence radius condition above is satisfied if $\int_{ x =0}^1\|\alpha(x)\|_{\mathfrak g} <2.4$ holds.
(Cf. \cite{L3}.)
In this viewpoint there is little geometry, and there is no actual straightening of the path $ x \mapsto g( x )$.
\snewpage

\textbf{The Maurer--Cartan viewpoint.}
Nevertheless, we can examine what happens if we try to do such a straightening.
Assume now that there is a smooth extension $\tilde{g}:[0,a]\times[0,1]\rightarrow G $ such that
$g( x )=\tilde{g}(0, x )$, while $\tilde{g}( t ,0)=1$ and $\tilde{g}( t ,1)=h$ throughout.
Beside
\[A( t , x )=\tilde{g}( t , x )^{-1}\frac{\partial}{\partial   x }\tilde{g}( t , x ),\]
we can also define
 \[B( t , x )=\tilde{g}( t , x )^{-1}\frac{\partial}{\partial  t }\tilde{g}( t , x ).\]
These quantities satisfy the (special case of the) Maurer--Cartan equation
\begin{equation}
\frac{\partial}{\partial  t }A( t , x )-\frac{\partial}{\partial  x }B( t , x )
-[A( t , x ),B( t , x )]
=0.
\plabel{MC1}
\end{equation}
(This is of Maurer's form. In a more abstract language,
$\eta ( t , x )=A( t , x )\mathrm d x +B( t , x )\mathrm d t $
satisfies
\[\mathrm d\eta ( t , x )+\frac12[\eta ( t , x ),\eta( t , x )]=0.\]
This is Cartan's form. Cf. Maurer \cite{Mau}, Cartan \cite{Car}, Bourbaki \cite{Bou}.)
The invariance conditions on the ends translate to
\begin{equation}
B( t ,0)=0 \qquad\text{and}\qquad B( t ,1)=0\qquad\text{for}\qquad t >0.
\plabel{MC2}
\end{equation}
In turn, if $A$ and $B$ are given so, then such  a $\tilde{g}$ can be constructed uniquely
(invariance for multiplication on the left is countered by the condition $\tilde{g}(0,0)=1$.)
The Maurer--Cartan equation can be imagined as a process, when there is a non-commutative mass
$A( t , x )$, whose change of rate in time $ t $ comes from the gradient of the (inverse) current $B( t , x )$
plus an interaction term as the current passes through the mass and changes it.

\textbf{The heat prescription.}
It may be reasonable to choose the current as the gradient of the mass itself, such that
\begin{equation}
B( t , x )=k\frac{\partial}{\partial  x }A( t , x ),
\plabel{DP}
\end{equation}
where $k>0$ is a diffusion parameter.
(Terminology: The multiplicative inverse $m=k^{-1}$ is the particular mass of the diffusion; ``heavier particles'' diffusing more slowly.)
Then the Maurer--Cartan equation yields
\begin{equation}
\frac{\partial}{\partial  t }A( t , x )-k\frac{\partial^2}{\partial  x ^2}A( t , x )
-\left[A( t , x ),k\frac{\partial}{\partial  x }A( t , x )\right]
=0.
\plabel{eq:hmc1}
\end{equation}
Beside the initial condition,
 \begin{equation}
A(0, x )=\alpha( x ) ,
\plabel{eq:hmc2}
\end{equation}
by the invariance of the ends, the boundary conditions
\[
\left.\frac{\partial}{\partial  x }A( t , x )\right|_{ x =0}=0
\quad\text{and}\quad
\left.\frac{\partial}{\partial  x }A( t , x )\right|_{ x =1}=0
\qquad\text{for}\qquad t >0\]
must be introduced.
In general, the Maurer--Cartan equation \eqref{MC1}--\eqref{MC2} is the natural consistency criterion for any straightening,
while the diffusion prescription \eqref{DP} is just one possible but relatively natural way to carry it out.

\snewpage

\textbf{The $t\rightarrow+\infty$ limit.}
It must be clear that the diffusion prescription \eqref{DP} is a natural not-to-think prescription,
 which is therefore natural and widely applicable.
One cannot hope to have a straightening in finite time (even the case of commutative Lie algebras shows that),
 but in appropriate circumstances for $ t \rightarrow+\infty$ (which supposes $a=+\infty$)
 homogenization  occurs in the limit.
More precisely one can argue as follows:

Having a norm ready, one can formulate conditions for this to be happen:
We suppose that $a=+\infty$.
Then one can define the initial [variation of] mass as
\[M_I(\alpha)=\int_{ x =0}^1\|\alpha (x) \| ,\]
and the [variation of] generated mass as
\[M_G(A)=\int_{( t , x )\in([0,+\infty)\times[0,1]}\|[A( t , x ),B( t , x ) ]\|\,\mathrm d t \,\mathrm d x .\]
Now, if
\begin{equation}
M_I(\alpha)+M_G(A)<+\infty,
\plabel{eq:sume}
\end{equation}
then,   $A_{ t }:[0,1]\rightarrow\mathfrak g$
(defined by $A_ t ( x )=A( t , x )$) gets increasingly homogenized as $ t \rightarrow+\infty$, and
\begin{equation}
H=\lim_{ t \rightarrow+\infty}\int_{ x =0}^1 A_{ t }( x ) \,\mathrm d x
\plabel{eq:limi}
\end{equation}
will have the property that $h=\exp H$.
One can also see that for $\tau\in[0,+\infty)$
\[\|H\|\leq \int_{ x =0}^1\|A_\tau( x )\|\,\mathrm d x  +
\int_{( t , x )\in(\tau,+\infty)\times[0,1]}\|[A( t , x ),B( t , x ) ]\|\,\mathrm d t \,\mathrm d x .\]
In particular, getting it to $\tau=0$, we find that
\[\|H\|\leq M_I(\alpha)+M_G(A).\]
(In this case an actual straightening of a path is exhibited by a evolution as $t\rightarrow+\infty$.)

Let us now consider what guarantees the existence of $A$ from having the initial $\alpha$, and such that \eqref{eq:sume} holds.
Regarding this, the formal approach will be of use again.

For example, we will demonstrate, in terms of Banach--Lie norms,
\begin{theorem}\plabel{thm:ex}
If $M_I(\alpha)\equiv\int \|\alpha \|_{\mathfrak g} <1$,
then $A$ as a solution of the non-commutative heat equation for $(t,s)\in[0,+\infty)\times [0,1]$ can be constructed, such that
\[M_I(\alpha)+M_G(A)\leq 2-2\sqrt{1-M_I(\alpha)}.\]
\end{theorem}
This statement will, in fact, be demonstrated in multiple versions, depending on the spacial domain.
This leads to
\snewpage

\textbf{Variations on the theme.}
Next we discuss the relaxation of some conditions of the situation above.

Firstly, we can replace the path parameter space $[0,1]$ by $[0,+\infty)$ or $\mathbb R\equiv(-\infty,+\infty)$.
Here we can assume $\alpha$ to be rapidly decreasing.
More or less, the same discussion applies.
In the latter cases $A_t$ cannot be expected to converge in distribution
(in case of a process of bounded mass), however, if the distribution is rescaled, then it does.
In fact, the technically easiest case is when the domain is $\mathbb R$ (as the heat propagation can be written down simply).

Secondly, although it introduces issues in the technical formulation, no smoothness (or rapid decrease) condition on $\alpha$ is necessary.
In fact, the density $\alpha( x )\,\mathrm d x $ can be replaced by a $\mathfrak g$ valued measure $\alpha$.

Thirdly, situation can be adapted to the case when we have a ``multicomponent system''
 which is roughly speaking is when $\mathfrak g$ is a positively graded vector space.
Then $\boldsymbol k:\mathfrak g\rightarrow\mathfrak g$ is the grading map
 which multiples by the (variable)   diffusion parameter grade-wise.
Its multiplicative inverse $\boldsymbol m=\boldsymbol k^{-1}:\mathfrak g\rightarrow\mathfrak g$
 multiples by the (variable) particular mass grade-wise.
Some natural choices are as follows.
Assume that $\mathfrak g$ is graded by $\mathbb N\setminus\{0\}$ (like the ambient Lie algebra of the formal noncommutative masses).
Let $\gr$ be the grading map (acting by multiplication with $\mathbb N\setminus\{0\}$ gradewise).
Then a natural choice is $\boldsymbol m=m^*\mathrm e^{\beta\cdot \gr}$ with $m^*\in(0,+\infty)$, $\beta\in\mathbb R$
(making  $\boldsymbol k=k^*\mathrm{e}^{-\beta\cdot\gr}$).
As $t\rightarrow +\infty$, for  $D=[0,+\infty)$ or $D=(-\infty,+\infty)$, we will not have a homogenization
in the limit but a kind of imperfect fractionalization by the ``particular mass'' $\boldsymbol m$.
This imperfect fractionalization is not suitable for any of our Lie expansions yet.
For this we have to take a secondary limit in limits.
The limit $\beta\rightarrow0$ diminishes the fractionalization, and we reobtain the Magnus expansion.
The limits $\beta\rightarrow+\infty$ and $\beta\rightarrow+\infty$ amplify the fractionalization, and we obtain the (symmetric) Wilcox expansions.
(In keeping track, formal solutions will be useful.)
More precisely, we have, for example, the following cases:

(i) If $D=[0,1]$, then, as $t\rightarrow+\infty$, homogenization occurs (under favorable circumstances)  yielding
\[h=\exp H.\]

(ii) If $D=[0,+\infty)$, then, as $t\rightarrow+\infty$, $\beta\rightarrow +\infty$, it yields (under favorable circumstances)
\[h=\ldots(\exp H_3)(\exp H_2)(\exp H_1),\]
where $H_i$ is of grade $i$.
Meanwhile $\beta\rightarrow -\infty$ yields (under favorable circumstances)
\[h=(\exp H_1)(\exp H_2)(\exp H_3)\ldots,\]
where $H_i$ is of grade $i$.
(It is not necessary to use this limit:
If the domain is $D=(-\infty,0]$, then the order of the components reverses.)

(iii) If  $D=(-\infty,+\infty)$, then, as $t\rightarrow+\infty$,
$\beta\rightarrow +\infty$, it yields (under favorable circumstances)
\[h= (\exp H_1/2)(\exp H_2/2(\exp H_3/2))\ldots(\exp H_3/2)(\exp H_2/2)(\exp H_1/2), \]
where $H_i$ is of grade $i$ again.
 If
$\beta\rightarrow -\infty$, then it yields (under favorable circumstances)
\[h=\ldots(\exp H_3/2)(\exp H_2/2)(\exp H_1/2)(\exp H_1/2)(\exp H_2/2)(\exp H_3/2)\ldots \]
where $H_i$ is of grade $i$ again.

In fact, in the latter cases $H_i$ is the grade $i$ part of $H$ as in \eqref{eq:limi} but integrated on the given domain,
 and then the limit $\beta\rightarrow \pm\infty$ taken.
If we apply these limits to the to formal noncommutative masses, then we obtain
 information regarding the Wilcox expansion and its symmetric versions.

\begin{remark}
Although the Maurer--Cartan equation is sensitive to the choice of left or right invariant vector fields taken,
on the Banach algebraic or formal level, if $k$ is constant,  the heat prescription yields the multiplicative heat equation
\[\frac{\partial}{\partial  t }\tilde g( t , x )=k \frac{\partial^2}{\partial  x ^2}\tilde g( t , x )
-k \frac{\partial }{\partial  x  }\tilde g( t , x ) \tilde g( t , x )^{-1}\frac{\partial }{\partial  x  }\tilde g( t , x ).\]
This is left-right (i. e. transposition) invariant.
This form is, however, not easily adaptable to the case with variable $k$.
In fact, as gradings are generally not conjugation (or $\ad$) invariant, transposition invariance gets broken.
\qedremark
\end{remark}

Fourthly, the domain $[0,1]$ with reflective boundary conditions can be replaced by the periodic boundary condition.
In general, the Mauer--Cartan condition without the invariance on the ends gives, for $\tau>0$,
in terms of time-ordered exponentials,
 \begin{multline*}
 \exp_{\mathrm R}( x \in[0,1]\mapsto A(\tau, x ) )=\\=
  \exp_{\mathrm R}( t \in[-\tau,0]\mapsto -B(- t ,0) )
 \exp_{\mathrm R}( x \in[0,1]\mapsto A(0, x ) )
 \exp_{\mathrm R}( t \in[0,\tau]\mapsto B( t ,1) ).
 \end{multline*}
In the case of periodic boundary conditions, $B( t ,0)=B( t ,1)$, we have conjugation here,
 \[
 \exp_{\mathrm R}( x \in[0,1]\mapsto A(\tau, x ) )=
  (F_\tau)^{-1}\cdot
 \exp_{\mathrm R}( x \in[0,1]\mapsto A(0, x ) )\cdot
 F_\tau.
 \]

In the case of the heat prescription,
 the cumulative mass of the ``boundary flux''  $t\mapsto B( t ,0)=B( t ,1)$ can be bounded by the cumulative total mass;
 thus in the case of finite mass \eqref{eq:sume}, in the infinite limit, this yields
\[h=F_\tau\cdot(\exp H)\cdot (F_\tau)^{-1}.\]
Consequently,
\[h= \exp \tilde H ,\]
where
\[\tilde H=F_\tau\cdot  H \cdot (F_\tau)^{-1}=\exp_{\mathrm R}( t \in[0,+\infty)\mapsto \ad B( t ,0)  ) \,\tilde H.\]
Applied to the formal non-commutative masses, the latter form also must give a presentation for the Magnus expansion.
\snewpage

Using these ideas one can prove
\begin{prop}\plabel{thm:magnustie}
\begin{equation}
\sum_{n=1}^\infty T^n\cdot\|\mu_{\mathrm R,n}(\alpha)\|_{\mathfrak g}
\stackrel{\forall T^n}{\leq} 2-2\sqrt{1-T\cdot\left(\int\|\alpha \|_{\mathfrak g} \right)}.
\plabel{eq:magnustie}
\end{equation}
In particular, the Magnus expansion is  convergent if $\int\|\alpha\|_{\mathfrak g}<1$.
\end{prop}
\noindent\textit{Note.} Using some intuitive arguments, one can obtain that  the Magnus expansion convergent if
\begin{equation}
\int\|\alpha\|_{\mathfrak g}<4-2\sqrt2=1.1715\ldots
\plabel{eq:magextra}
\end{equation}
holds; and, in fact, we can do even better.
\qedremark

\begin{prop}\plabel{thm:zasstie}
\begin{equation}
\sum_{n=1}^\infty T^n\cdot\|\zeta^{\leftarrow}_{\mathrm R,n}(\alpha)\|_{\mathfrak g}
\stackrel{\forall T^n}{\leq}
 1- \sqrt{\left(2-T\cdot\left(\int\|\alpha \|_{\mathfrak g}\right)\right)^2-1} .
 \plabel{eq:zasstie}
\end{equation}
In particular, the   Wilcox expansion is convergent if $\int\|\alpha\|_{\mathfrak g}<2-\sqrt2=0.5857\ldots$ holds.

Similar statement holds with respect to the Wilcox expansion variant $\zeta^{\rightarrow}$.
\end{prop}
\noindent\textit{Note.} Using some intuitive arguments, one can obtain
\begin{equation}
 \sum_{n=1}^\infty T^n\cdot\|\zeta^{\rightarrow}_{\mathrm R,n}(\alpha)\|_{\mathfrak g}
\stackrel{\forall T^n}{\leq}
2- T\cdot\left(\int\|\alpha \|_{\mathfrak g} \right)-
  \sqrt{ 3\left(\frac43-T\cdot\left(\int\|\alpha \|_{\mathfrak g} \right)^2\right)-\frac43}.
 \plabel{eq:zasstieextra}
\end{equation}
Similar statement applies to the $\zeta^{\leftarrow}$ version.
According to this, the left  and right expanding Wilcox expansions are convergent if $\int\|\alpha\|_{\mathfrak g}<\frac23$ holds.
\qedremark

\begin{theorem}\plabel{thm:bowtie}
\begin{equation}
\sum_{n=1}^\infty T^n\cdot\|\eta^{\bowtie}_{\mathrm R,n}(\alpha)\|_{\mathfrak g}
\stackrel{\forall T^n}{\leq} 2-2\sqrt{1-T\cdot\left(\int\|\alpha \|_{\mathfrak g} \right)}.
\plabel{eq:bowtie}
\end{equation}

In particular,  the inward expanding symmetric Wilcox expansion is convergent 
if  $\int\|\alpha\|_{\mathfrak g}<1$ holds.

Similar statement holds for  the outward expanding symmetric  Wilcox expansion $\eta^{\leftrightarrow}$.
\end{theorem}
\noindent\textit{Note.} Using some intuitive arguments, one can obtain
\begin{equation}
\sum_{n=1}^\infty T^n\cdot\|\eta^{\leftrightarrow}_{\mathrm R,n}(\alpha)\|_{\mathfrak g}
\stackrel{\forall T^n}{\leq} 4- T\cdot\left(\int\|\alpha \|_{\mathfrak g} \right)-
  \sqrt{ 2\left(4-T\cdot\left(\int\|\alpha \|_{\mathfrak g} \right)\right)^2-16}.
\plabel{eq:bowtieextra}
\end{equation}
According to this, the  outward expanding symmetric Wilcox expansion is convergent if
$\int\|\alpha\|_{\mathfrak g}<4-2\sqrt2=1.1715\ldots$ holds.

[This is not claimed for $\eta^{\bowtie}$, but one would expect a similar result.]
\qedremark

\snewpage

\section{The formal solution to the non-commutative heat equation}\plabel{sec:formheat}
\textbf{The construction of the formal solutions.}
The formal solution is a solution of the non-commutative equation where the
initial condition is replaced by
\[A^{(T)}(0, x )=T\cdot \alpha( x ) \]
where $T$ is a formal variable, supposedly with the same initial and boundary value conditions.
(Here the part `$(T)$' of the notation is not strictly required, and could be omitted, but, in what follows, we will use the
indicators `$(T)$' liberally just to indicate that we have something dealing with the formal approach.)
Here $\alpha$ is understood to be a Banach or Banach--Lie algebra smooth valued function.
This latter condition can be changed to being a Lebesgue--Bochner integrable function,
in particular, of $L^1$ in variation norm with respect to a Banach algebraic norm Banach--Lie algebraic norm.
In fact, $\alpha(x)\, \mathrm dx$ can be replaced by a  Banach or Banach--Lie algebra valued (interval) measure of finite variation.
Actually, the formal solutions which we consider will be directly constructed.
For that reason, with respect to the choice of spacial domains
$D=[0,1]$, $[0,+\infty)$, $\mathbb R$, $\mathbb R/\mathbb Z$, we will consider
the heat propagation functions (heat kernels)
\[K_{[0,1] }(x,y,t;k)=\sum_{p\in2\mathbb Z}\frac{1}{2\sqrt{\pi k t}}\exp\left(-\frac{(x-y-p)^2}{4kt}\right)
+\sum_{p\in2\mathbb Z}\frac{1}{2\sqrt{\pi k t}}\exp\left(-\frac{(x+y+p)^2}{4kt}\right);
\]
\[K_{[0,+\infty)}(x,y,t;k)=\frac{1}{2\sqrt{\pi k t}}
  \exp\left(-\frac{(x-y)^2}{4kt}\right)+ \frac{1}{2\sqrt{\pi k t}}\exp\left(-\frac{(x+y)^2}{4kt}\right);  \]
\[K_{\mathbb R}(x,y,t;k)=\frac{1}{2\sqrt{\pi k t}}\exp\left(-\frac{(x-y)^2}{4kt}\right);\]
\[K_{ \mathbb R/\mathbb Z}(x,y,t;k)=\sum_{p\in \mathbb Z}\frac{1}{2\sqrt{\pi k t}}\exp\left(-\frac{(x-y-p)^2}{4kt}\right).\]
Here $x,y\in D$, $t>0$, $m=k^{-1}>0$.
\begin{commentx}

Remark: The probability density function of the standard normal distribution with mean $\mu$ and variance $\sigma$ is given by
\[N(x;\mu,\sigma)=\frac1{\sqrt{2\pi\sigma^2}}\exp\left(-\frac{(x-\mu)}{2\sigma^2}\right).\]
\end{commentx}
Then the formal solution we will consider is of shape
\[A^{(T)}(t,x)=\sum_{n=1}^\infty T^n\cdot A_n(t,x),\]
where the solution basically propagates along the heat kernels but it gets interacted to higher powers of $T$.
More concretely, for $t>0$:
\begin{equation}
A_1(t,x)= \int_{y\in D} K_D(x,y,t;k)\alpha(y).
\plabel{eq:ren1}
\end{equation}
If $A_1(t,x),\ldots,A_{n-1}(t,x)$ are constructed, then the generated mass of $n$th order is
\begin{equation}
A_{n}^{\gen}(t,x)=\sum_{j=1}^{n-1}\left[A_j( t , x ),k\frac{\partial}{\partial  x }A_{n-j}( t , x )\right].
\plabel{eq:ren2}
\end{equation}
Then the propagating mass of $n$th order is
\begin{equation}
A_{n} (t,x)= \int_{s\in(0,t)}\int_{y\in D} K_D(x,y,t-s;k) A_{n}^{\gen}(s,y)\,\mathrm dy\,\mathrm ds.
\plabel{eq:ren3}
\end{equation}
More precisely, this is the situation in the case when $k>0$ is constant.

However, as we have indicated, we can also consider a more general case, when
 our algebra is positively graded as a vector space but not necessarily as a (Lie) algebra.
Here we spell out this case.
For the sake of simplicity, we consider only the case when the grading take values in a countable positive set
 $\mathcal P(\boldsymbol m)$, but situation can be adapted continuous positive gradings.
We will use the notation $\proj_m^{\boldsymbol m}$ as the projection to grade $m$.
We will also assume that the norm is gradewise induced, i. e.
 \[\|X\|=\sum_{m\in\mathcal P(\boldsymbol m) }\|\proj^{\boldsymbol m}_mX\|.\]
(This is sufficient with respect to the applications to the Lie series, as we want to use formal noncommutative masses
 $\alpha$ with respect to the natural ``polynomial'' grading $\gr$ and prescriptions $\boldsymbol m=m^*\mathrm e^{\beta\cdot \gr}$.)
Then \eqref{eq:ren1}, \eqref{eq:ren2}, \eqref{eq:ren3} must be replaced by their counterparts
 \begin{equation}
A_1(t,x)= \int_{y\in D} K_D(x,y,t;\tfrac1m)\sum_{m\in\mathcal P(\boldsymbol m) }\proj^{\boldsymbol m}_m\alpha(y),
\plabel{eq:ren11}
\end{equation}
\begin{equation}
A_{n}^{\gen}(t,x)=\sum_{j=1}^{n-1}\left[A_j( t , x ),
\tfrac1m\sum_{m\in\mathcal P(\boldsymbol m) }\proj^{\boldsymbol m}_m\frac{\partial}{\partial  x }A_{n-j}( t , x )\right],
\plabel{eq:ren22}
\end{equation}
\begin{equation}
A_{n} (t,x)= \int_{s\in(0,t)}\int_{y\in D} K_D(x,y,t-s;\tfrac1m)
\tfrac1m\sum_{m\in\mathcal P(\boldsymbol m) }\proj^{\boldsymbol m}_m
A_{n}^{\gen}(s,y)\,\mathrm dy\,\mathrm ds.
\plabel{eq:ren33}
\end{equation}

In order to obtain control over this construction,
 we will consider the solution as a superposition generated by appropriately placed Dirac delta functions.
\snewpage

\textbf{Analytic control  in the case $D=\mathbb R$.}
Let us first consider the case when $D=\mathbb R$.
Let us consider the case when $T\cdot\alpha=T\cdot Y_1\cdot \boldsymbol\delta_{y_1}+T\cdot Y_2\cdot \boldsymbol\delta_{y_2}$, $y_1<y_2$,
 where $Y_1$ and $Y_2$ are Lie algebra elements of (vector space) grade $m_1$ and $m_2$ respectively.
The initial mass propagates  according to
\[A_1(t,x)= K_{\mathbb R}(x,y_1,t;\tfrac1{m_1})\cdot Y_1
+ K_{\mathbb R}(x,y_2,t;\tfrac1{m_2})\cdot Y_2.\]
The primary interaction term, which is the mass directly generated from the initial mass is according to
\begin{multline*}
A_2^{\gen}(t,x)= K_{\mathbb R}(x,\tfrac{m_1y_1+m_2y_2}{m_1+m_2},t;\tfrac1{m_1+m_2})\\
\frac14\,{\frac {\sqrt {m_1 m_2}  }{\sqrt {\pi }{t}^{3/2}\sqrt {m_1+{m_2}}}{\exp{\left(-\frac14\,{
\frac {m_2\,m_1\, \left( { y_2}-{ y_1} \right) ^{2}}{t\,
 \left( m_1+{m_2} \right) }}\right)}}}\cdot
\left(  y_2- y_1 \right)[Y_1,Y_2].
\end{multline*}
This directly generated mass density (apart from $T^2$) is positively proportional to $[Y_1,Y_2]$.

\begin{multline*}
\int_{t\in(0,+\infty)}\int_{x\in\mathbb R}  A_2^{\gen}(t,x)\,\mathrm dx\,\mathrm dt=\\
=
\int_{t\in(0,+\infty)}
\frac14\,{\frac {\sqrt {m_1 m_2}  }{\sqrt {\pi }{t}^{3/2}\sqrt {m_1+{m_2}}}{\exp{\left(-\frac14\,{
\frac {m_2\,m_1\, \left( { y_2}-{ y_1} \right) ^{2}}{t\,
 \left( m_1+{m_2} \right) }}\right)}}}\cdot
\left(  y_2- y_1 \right)[Y_1,Y_2]
\,\mathrm dt\\
=
\left[ -\frac12\erf\left(\frac12\frac{(y_2-y_1)\sqrt {m_1 m_2}
}{\sqrt{t}\sqrt{m_1+m_2}}\right)\cdot [Y_1,Y_2]\right]^{+\infty}_{t=0}= \frac12[Y_1,Y_2].
\end{multline*}
(Here $\erf(\tau)=\frac{2}{\sqrt\pi}\int_{t=0}^\tau\exp(-t^2)\,\mathrm dt$.)
Another way to see the same thing is via
\begin{multline*}
\int_{x\in\mathbb R}  \int_{t\in(0,+\infty)}A_2^{\gen}(t,x)\,\mathrm dt\,\mathrm dx
=
\int_{x\in\mathbb R}
\frac{(y_2-y_1)\sqrt{m_1m_2}}{2\pi\left(m_1(x-y_1)^2+m_2(x-y_2)^2\right)}
[Y_1,Y_2]
\,\mathrm dt\\
=
\left[ \frac1{2\pi}\arctan\left(\frac{m_1(x-y_1)+m_2(x-y_2)}{(y_2-y_1)\sqrt{m_1m_2}}\right)  \cdot [Y_1,Y_2]\right]^{+\infty}_{t=-\infty}= \frac12[Y_1,Y_2].
\end{multline*}

Similar argument applies for mass generated from generated masses
$T^{j}\cdot Y_1\cdot \boldsymbol\delta_{(s,y_1)}$ and $T^{n-j}\cdot Y_2\cdot \boldsymbol\delta_{(s,y_2)}$, $y_1<y_2$,
placed at the same time $s$.
Ultimately, the mass $T^n\cdot\frac12[Y_1,Y_2]$ will be distributed by a probability measure as generated mass, onward from time $s$.
When we consider the mass generated from
$T^{j}\cdot Y_1\cdot \boldsymbol\delta_{(s_1,y_1)}$ and $T^{n-j}\cdot Y_2\cdot \boldsymbol\delta_{(s_2,y_2)}$,
 such that $s_1\neq s_2$ (different times), then we can argue as follows.
If $s_1<s_2$, the first we wait until $T^{j}\cdot Y_1\cdot \boldsymbol\delta_{(s_1,y_1)}$ propagates to time $s_2$,
 then we disintegrate it to delta functions $T^{j}\cdot Y_1\cdot \boldsymbol\delta_{(s_1,\tilde y_1)}$ , and consider
 the interaction terms with  $T^{n-j}\cdot Y_2\cdot \boldsymbol\delta_{(s_2,y_2)}$.
After the propagation, both $\tilde y_1<y_2$ and $\tilde y_1>y_2$ can occur (with complementary nonzero probabilities).
Ultimately $T^n\cdot\frac12[Y_1,Y_2]$ and $T^n\cdot\frac12[Y_2,Y_1]$ will be distributed with complementary probabilities;
 actually none of them with zero density, thus cancellations between them occur them even as they are created.
Thus in variation, the generated mass is strictly less than $T^n\cdot\frac12\|[Y_1,Y_2]\|$ or zero.
(Here $\|\cdot\|$ may mean Banach algebra norm of Banach-Lie algebra norm either.)
Ultimately, we find that the mass generated from
 $T^{j}\cdot Y_1\cdot \boldsymbol\delta_{(s_1,y_1)}$ and $T^{n-j}\cdot Y_2\cdot \boldsymbol\delta_{(s_2,y_2)}$,
 independently from their placement, can be majorized in variation by $T^n\cdot\frac12\|[Y_1,Y_2]\|$.
\snewpage

The following discussion will specify to the case when $\|\cdot\|=\|\cdot\|_{\mathfrak g}$ is a Banach--Lie norm.
Let we define the overall norm majorizing function as
\[f(T)=T\cdot\left(\int\|\alpha\|_{\mathfrak g}\right)+\sum_{n=2}^\infty T^n\cdot
\left(\int_{t\in(0,+\infty)}\int \|A_n^{\gen}(t,x)\|_{\mathfrak g}\,\mathrm dx\,\mathrm dt\right).\]
Then the preceding discussion about the generated masses implies that
\begin{equation}
f(T) \stackrel{\forall T^n}{\leq} T\cdot\left(\int\|\alpha \|_{\mathfrak g}\right)+ \frac12\cdot\frac12\,f(T)^2
\plabel{eq:normest0}
\end{equation}
(A factor $\frac12$ appears at the commutator, and an other  factor $\frac12$ appears due to the spacial ordering in accounting.)
Then $f(T)$ gets majorized by the formal series $g(T)$, which is the solution of $g(T)=0$ and
\[g(T)=T\cdot\left(\int\|\alpha\|_{\mathfrak g}\right)+\frac12\cdot\frac12\, g(T)^2.\]
In concrete terms,
\begin{equation}
f(T) \stackrel{\forall T^n}{\leq} 2-2\sqrt{1-T\cdot\left(\int\|\alpha\|_{\mathfrak g}\right)}.
\plabel{eq:normest1}
\end{equation}
Putting $T=1$ here, we obtain the corresponding version of Theorem \ref{thm:ex}.

For $t>0$, we can consider the function $A_t^{(T)}$ given by $A^{(T)}_t(x)=A^{(T)}(t,x)$.
By (the extension of) the Maurer--Cartan equation, it holds that
\[\exp_{\mathrm R} (T\cdot\alpha)=\exp_{\mathrm R}(A^{(T)}_t).\]
We can also note that   $H(T)=\lim_{t\rightarrow +\infty}\int\,A^{(T)}_t(x)\mathrm dx$ exists, and, in fact,
\[H(T)=T\cdot\left(\int \alpha  \right)+\sum_{n=2}^\infty T^n\cdot
\left(\int_{t\in(0,+\infty)}\int  A_n^{\gen}(t,x) \,\mathrm dx\,\mathrm dt\right).\]
Moreover, if we rescale $A^{(T)}_t$ by $\sqrt t$, then we find, as $t\rightarrow\infty$,
\begin{equation}
A^{(T)}_t\xrightarrow{\mathrm{rescaled}}A^{(T)}_\infty=\sum_{m\in\mathcal P} (\proj^{\boldsymbol m}_m H)\cdot K_{\mathbb R}(x,0,1;\tfrac1m).
\plabel{eq:resc}
\end{equation}
It will hold that
\[\exp_{\mathrm R} (T\cdot\alpha)=\exp_{\mathrm R}(A^{(T)}_\infty).\]

\snewpage

\textbf{The $\beta=-\infty,0,+\infty$ cases for $D=\mathbb R$.}
Now we assume that we deal with formal noncommutative masses as measures and $\boldsymbol m=m^*\mathrm e^{\beta\cdot \gr}$.

Let us consider now the special case when $k$ is constant ($\beta=0$).
With $ A^{(T)}_\infty$ is just $H(T)$ distributed by a probability measure,
\[\exp_{\mathrm R} (T\cdot\alpha)=\exp H(T).\]
As we have discussed, in this formal case $H(T)$ must realize  the (formal) Magnus expansion,
\[H(T)=\sum_{n=1}^\infty T^n\cdot\mu_{\mathrm R,n}(\alpha).\]
From \eqref{eq:normest1},
we obtain
Proposition \ref{thm:magnustie}.
This is weaker than \eqref{eq:magmo} but it can be counted as a natural ``trivial estimate''
in the sense that it ignores the finer, actual, combinatorial details of the $\mu_n$.
Beside the actual convergence of the Magnus series, it also says that if a noncommutative mass of cumulative norm less than $1$ is left to be
 diffused multiplicatively in a uniform manner, then it homogenizes to the mass of the Magnus expansion.

Next we consider the case $\beta\rightarrow+\infty$.
(For this reason all previous expressions are imagined to be endowed by the indicator `$[\beta]$'.)
As it was indicated as before, in the limit the distribution of $A^{(T)[\beta]}_t$ fractionalizes by $\gr$;
we can (cut off higher powers of $T$ and) reparametrize the spacial variable on the RHS of \eqref{eq:resc}.
Then
\[H(T)=\lim_{\beta\rightarrow+\infty}H^{[\beta]}(T)=\lim_{\beta\rightarrow+\infty}\int A_\infty^{(T)[\beta]}
\]
exits, and the uniform estimates of \eqref{eq:normest1} are preserved; and
\begin{multline*}
\exp_{\mathrm R} (T\cdot\alpha)=\exp(\proj^{\gr}_1 H(T)/2) \exp(\proj^{\gr}_2 H(T)/2)  \exp(\proj^{\gr}_3 H(T)/2) \ldots
\\\ldots\exp(\proj^{\gr}_3 H(T)/2) \exp(\proj^{\gr}_2 H(T)/2)  \exp(\proj^{\gr}_1 H(T)/2)
\end{multline*}
holds.
Arguing in the same manner as before, we obtain
\[H(T)=\sum_{n=1}^\infty T^n\cdot\eta^{\bowtie}_{\mathrm R,n}(\alpha)\]
and the $\eta^{\bowtie}$-part of  Theorem \ref{thm:bowtie}.

The limit $\beta\rightarrow-\infty$ serves $\eta^{\leftrightarrow}_{\mathrm R,n}(\alpha)$ analogously;
leading to the $\eta^{\leftrightarrow}$-part of  Theorem \ref{thm:bowtie}.
\snewpage

\medskip
\noindent\textit{Sketch of improvement.}
This applies in the formal case when $\beta\rightarrow-\infty$.
Then the mass of higher grade diffuses in much higher rate than mass of lower grade.
Let us consider mass generation occuring from mass of grade $n'$ and $n$ with $2\leq n'\leq n$.
Then on the scale of grade level $n$ diffusion, the generated mass of grade $n$ or less
 gets introduced in the diffusing mass  of grade $n$ very slowly, at a sort of glacial rate.
This means that generated mass gets introduced into the process as $\boldsymbol \delta_0(x) Y'$
 compared to the diffused mass $K_{\mathbb R}(x,0,t;m) Y$.
In this case the generated mass
 is $\frac18[Y',Y]$ and $\frac18[Y,Y']$ distributed by probability measures; altogether of bounded by
 total variation $\frac14\|[Y,Y']\|$.
This means that the limit
\[f(t)=\lim_{\beta\rightarrow -\infty} f^{[\beta]}(T)\]
will be majorized by the solution $g(T)$ of $g(T)=0$ and
\begin{multline*}
g(T)=T\cdot\left(\int\|\alpha\|_{\mathfrak g}\right)+\frac12\,T^2\cdot\frac12\left(\int\|\alpha\|_{\mathfrak g}\right)^2 +\\
+\frac12\,T\cdot\left(\int\|\alpha\|_{\mathfrak g}\right)\left( g(T)-T\cdot\left(\int\|\alpha\|_{\mathfrak g}\right)\right)
+\frac14\cdot\frac12\left( g(T)-T\cdot\left(\int\|\alpha\|_{\mathfrak g}\right)\right)^2.
\end{multline*}
In this way, we obtain  \eqref{eq:bowtieextra}.
\qedremark

\snewpage

\textbf{The case $D=[0,+\infty)$.}
Having already seen a similar discussion, we see that it sufficient to treat the mass generation
coming from delta functions placed at equal times.
Again, let us consider the case when $T\cdot\alpha=T\cdot Y_1\cdot \boldsymbol\delta_{y_1}+T\cdot Y_2\cdot \boldsymbol\delta_{y_2}$, $y_1<y_2$,
 where $Y_1$ and $Y_2$ are Lie algebra elements of (vector space) grade $m_1$ and $m_2$ respectively.
The primary interaction term, which is the mass directly generated from the initial mass is according to
\begin{multline*}
A_2^{\gen}(t,x)= \left(K_{\mathbb R}(x,\tfrac{m_1y_1+m_2y_2}{m_1+m_2},t;\tfrac1{m_1+m_2})-
K_{\mathbb R}(x,\tfrac{-m_1y_1-m_2y_2}{m_1+m_2},t;\tfrac1{m_1+m_2})\right)
\\
\frac14\,{\frac {\sqrt {m_1 m_2}  }{\sqrt {\pi }{t}^{3/2}\sqrt {m_1+{m_2}}}{\exp{\left(-\frac14\,{
\frac {m_2\,m_1\, \left( { y_2}-{ y_1} \right) ^{2}}{t\,
 \left( m_1+{m_2} \right) }}\right)}}}\cdot
\left(  y_2- y_1 \right)[Y_1,Y_2]\\
+ \left(K_{\mathbb R}(x,\tfrac{-m_1y_1+m_2y_2}{m_1+m_2},t;\tfrac1{m_1+m_2})-
K_{\mathbb R}(x,\tfrac{m_1y_1-m_2y_2}{m_1+m_2},t;\tfrac1{m_1+m_2})\right)
\\
\frac14\,{\frac {\sqrt {m_1 m_2}  }{\sqrt {\pi }{t}^{3/2}\sqrt {m_1+{m_2}}}{\exp{\left(-\frac14\,{
\frac {m_2\,m_1\, \left( { y_2}+{ y_1} \right) ^{2}}{t\,
 \left( m_1+{m_2} \right) }}\right)}}}\cdot
\left(  y_2+ y_1 \right)[Y_1,Y_2].
\end{multline*}

Some computation yields that the first summand distributes the mass
\[\frac1\pi\arctan\left(\frac{m_2y_2+m_1y_1}{\sqrt{m_1m_2}(y_2-y_1)}\right)[Y_1,Y_2]\]
according to a probability distribution;
while the second summand distributes the mass
\[\frac1\pi\arctan\left(\frac{m_2y_2-m_1y_1}{\sqrt{m_1m_2}(y_2+y_1)}\right)[Y_1,Y_2].\]
(There are no sign changes in the integrands!)
If $m_1=m_2$, then the mass $\frac12[Y_1,Y_2]$ is distributed.
 This leads to the variant of Theorem \ref{thm:ex} for $D=[0,+\infty)$ by the arguments we have used before.
In general, however, we cannot do much better than to have the estimate
\[\int_{t\in(0,+\infty)}\int_{x\in[0,+\infty)}  \|A_2^{\gen}(t,x)\|\,\mathrm dx\,\mathrm dt\leq \|[Y_1,Y_2]\|.\]
Ultimately, by similar arguments as before, we find that the mass generated from
 $T^{j}\cdot Y_1\cdot \boldsymbol\delta_{(s_1,y_1)}$ and $T^{n-j}\cdot Y_2\cdot \boldsymbol\delta_{(s_2,y_2)}$,
 independently from their placement, can be majorized in variation by $T^n\cdot\frac12\|[Y_1,Y_2]\|$
 if $m_1=m_2$, and it can be majorized by $T^n\cdot \|[Y_1,Y_2]\|$ in general.
Thus, in the fixed $k=m^{-1}$ case, we will a similar estimate for the Magnus expansion as in \eqref{eq:magnustie}.
In the general case, we only will have, regarding the generated masses, that
\[f(T) \stackrel{\forall T^n}{\leq} T\cdot\left(\int\|\alpha\|_{\mathfrak g}\right)+ \frac12f(T)^2\]
This leads to
\begin{equation*}
f(T) \stackrel{\forall T^n}{\leq} 1- \sqrt{1-T\cdot2\left(\int\|\alpha\|_{\mathfrak g}\right)}.
\end{equation*}

When we apply this argument to the Wilcox expansion, however, we can do better.
It is sufficient to note only that the mass of grade $2$ gets created not in variation $\|[Y_1,Y_2]\|$ but in
 variation $\frac12\|[Y_1,Y_2]\|.$
Therefore, $f(T)$ gets majorized by the solution $g(T)$ of $g(T)=0$ and
\[g(T)=T\cdot\left(\int\|\alpha\|_{\mathfrak g}\right)-\frac14\cdot T^2\cdot\left(\int\|\alpha\|_{\mathfrak g}\right)^2 +\frac12\, g(T)^2.\]

Using the same arguments as before,
the secondary limit $\beta\rightarrow+\infty$
yields the $\zeta^{\leftarrow}$-part of Theorem \ref{thm:zasstie}.

The limit $\beta\rightarrow-\infty$ deals with $\zeta^{\rightarrow}_{\mathrm R,n}(\alpha)$ analogously,
 yielding the $\zeta^{\rightarrow}$-part of Theorem \ref{thm:zasstie}.

\medskip
\noindent\textit{Sketch of improvement.}
This, again, applies in the formal case when $\beta\rightarrow-\infty$.
We can use our arguments from before in the case of mass created from grades $2\leq n'\leq n$.
This means that the generated mass gets introduced into the process as $\boldsymbol \delta_0(x) Y'$
 compared to the diffused mass $K_{[0,+\infty)}(x,0,t;m) Y$.
In this case the mass generated from this
 is $\frac12[Y',Y]$  distributed by a probability measure; altogether of bounded by
 total variation $\frac12\|[Y,Y']\|$.
This means that the limit
\[f(t)=\lim_{\beta\rightarrow -\infty} f^{[\beta]}(T)\]
will be majorized by the solution $g(T)$ of $g(T)=0$ and
\begin{multline*}
g(T)=T\cdot\left(\int\|\alpha\|_{\mathfrak g}\right)+\frac12\,T^2\cdot\frac12\left(\int\|\alpha\|_{\mathfrak g}\right)^2 +\\
+ T\cdot\left(\int\|\alpha\|_{\mathfrak g}\right)\left( g(T)-T\cdot\left(\int\|\alpha\|_{\mathfrak g}\right)\right)
+\frac12\cdot\frac12\left( g(T)-T\cdot\left(\int\|\alpha\|_{\mathfrak g}\right)\right)^2.
\end{multline*}
In this way, we obtain  \eqref{eq:zasstieextra}.
But, the machineries of the right- are left-expanding Wilcox series are equivalent as passing to
the (negative) transpose measure shows.
\qedremark

\snewpage

\textbf{The case $D=[0,1]$.}
Here we will restrict to the case of constant diffusion rate.
Due to symmetry reasons, it is sufficient only the case
   $T\cdot\alpha=T\cdot Y_1\cdot \boldsymbol\delta_{y_1}+T\cdot Y_2\cdot \boldsymbol\delta_{y_2}$,
 with $0\leq y_1<y_2\leq1$, where $Y_1$ and $Y_2$ are Lie algebra elements.
After some rearrangements, we find that the primarily generated mass is
\begin{align*}
A_2^{\gen}(t,x)=\sum_{s\in \mathbb Z}\sum_{r\in\mathbb Z}&
 \left( K_{\mathbb R}(x,\tfrac{y_1+y_2}{2}+s+2r,t;\tfrac1{2m})- K_{\mathbb R}(x,-\tfrac{y_1+y_2}{2}-s-2r,t;\tfrac1{2m}) \right)\cdot\\
&\cdot\frac14\,{\frac {\sqrt m }{\sqrt {\pi }{t}^{3/2}\sqrt {2}}{\exp{\left(-\frac14\,{
  \frac {\,m \left( { y_2}-{ y_1}+2s \right) ^{2}}{
  2 t\, }}\right)}}}\cdot
\left(  y_2- y_1+2s \right)[Y_1,Y_2]
\\&+ \left( K_{\mathbb R}(x,\tfrac{-y_1+y_2}{2}+s+2r,t;\tfrac1{2m})- K_{\mathbb R}(x, \tfrac{ y_1-y_2}{2}-s-2r,t;\tfrac1{2m}) \right)\cdot\\
&\cdot\frac14\,{\frac {\sqrt m }{\sqrt {\pi }{t}^{3/2}\sqrt {2}}{\exp{\left(-\frac14\,{
\frac {\,m \left( { y_2}+{ y_1}+2s \right) ^{2}}{
  2 t\, }}\right)}}}\cdot
\left(  y_2+ y_1+2s \right)[Y_1,Y_2].
\end{align*}
(Integrality has counted in the rearrangement.)
Here, for $x\in[0,1]$, we see that every summand is a nonnegative multiple of $[Y_1,Y_2]$.
We also know
\begin{equation*}
\int_{t\in(0,+\infty)}\int_{x\in[0,1]}  A_2^{\gen}(t,x)\,\mathrm dx\,\mathrm dt
= \frac12[Y_1,Y_2].
\end{equation*}
We know this not necessarily by a careful evaluation of the integral, but from the fact that the formal
Magnus expansion must work out in the second order.
(We could have also used this argument  in the constant diffusion rate cases  before.)
Therefore, in the mass generation $\frac12[Y_1,Y_2]$ is distributed by a probability measure.
We have the same rates as in \eqref{eq:normest1}, leading to \eqref{eq:magnustie}; or putting $T=1$ there, to the variant of
Theorem \ref{thm:ex}.
\snewpage

\textbf{The case $D=\mathbb R/\mathbb Z$.}
Let us consider the case of constant diffusion rate.
Before embarking on the actual computation, let us discuss the involved kernel in a bit more detail.
On physical grounds, it is quite obvious that   $K_{ \mathbb R/\mathbb Z}(x,y,t;k)$ is monotone in $\cos 2\pi(x-y)$.
There are intuitive arguments for this, but they are a bit painful to write down.
However, it is well-known that  $K_{ \mathbb R/\mathbb Z} $ can be rewritten by Jacobi's theta function as follows:
\[K_{ \mathbb R/\mathbb Z}(x,y,t;k)=\boldsymbol\vartheta_3\left(\pi(x-y),\exp(-4\pi^2k^2t^2 )\right),\]
where
\[\boldsymbol\vartheta_3(z,q)=1+2\sum_{n=1}^\infty q^{n^2}\cos(2nz)\]
(following the notation of Abramowicz, Stegun \cite{AS}, 16.27.3.);
see, e. g.,  Mardia, Jupp \cite{MJ}.
By Jacobi's  triple product formula
\[\sum_{n=-\infty}^\infty w^{2n}q^{n^2}=\prod_{m=1}^\infty(1-q^{2m})(1+w^2q^{2m-1})(1+w^{-2}q^{2m-1})\]
(see Whittaker,  Watson \cite{WW} and references therein), we see that
\[\boldsymbol\vartheta_3(z,q)=\prod_{m=1}^\infty(1-q^{2m})(1+ q^{4m-2}+2 q^{2m-1}\cos 2z). \]
From this, the monotonicity of $K_{ \mathbb R/\mathbb Z}(x,y,t;k)$ in $\cos 2\pi(x-y)$ is transparent.

Now, let us consider mass generation from
 $T\cdot\alpha=T\cdot Y_1\cdot \boldsymbol\delta_{y_1}+T\cdot Y_2\cdot \boldsymbol\delta_{y_2}$.
By translation invariance, it is actually sufficient to consider the case $y_1=-y$, $y_2=y$, $0<y<\frac14$.
The primarily   generated  mass is according to
\begin{multline*}
A_2^{\gen}(t,x)=\\= \sum_{n_1,n_2\in\mathbb Z}K_{\mathbb R}(x,\tfrac{ n_1 +n_2}{2},t;\tfrac1{2m })
\frac14\,{\frac {\sqrt {m }\cdot
\left(  2y+n_2-n_1 \right)  }{\sqrt {2\pi }{t}^{3/2} }{\exp{\left(-\frac14\,{
\frac {  m\left( { 2y } +n_2-n_1 \right) ^{2}}{2t
  }}\right)}}}[Y_1,Y_2]
\end{multline*}
\begin{multline*}
 = \underbrace{\left(\sum_{r\in2\mathbb Z}K_{\mathbb R}(x,\frac r2,t;\tfrac1{2m })\right)}_{K_{\mathbb R/\mathbb Z}(x,0,t;\frac1{2m})=}\cdot
\underbrace{\left(\sum_{n\in  2\mathbb Z}\frac14\,{\frac {\sqrt {m } \cdot
\left(  2y+n  \right) }{\sqrt {2\pi }{t}^{3/2} }{\exp{\left(-\frac14\,{
\frac {  m\left( { 2y } +n  \right) ^{2}}{2t
  }}\right)}}}\right)}_{S_{+}( y,0,t;\frac1{2m}):=}[Y_1,Y_2]
\end{multline*}
\begin{multline*}
 +\underbrace{\left(\sum_{r\in2\mathbb Z+1}K_{\mathbb R}(x,\frac r2,t;\tfrac1{2m })\right)}_{
 K_{\mathbb R/\mathbb Z}(x,\frac12,t;\frac1{2m})=}\cdot
\underbrace{\left(\sum_{n\in  2\mathbb Z}\frac14\,{\frac {\sqrt {m } \cdot
\left(  2y+n  \right) }{\sqrt {2\pi }{t}^{3/2} }{\exp{\left(-\frac14\,{
\frac {  m\left( { 2y } +n  \right) ^{2}}{2t
  }}\right)}}}\right)}_{S_{-}( y,0,t;\frac1{2m}):=}[Y_1,Y_2].
\end{multline*}
Then we see that $S_{+}( y,0,t;\frac1{2m})>0$ and $S_{-}(y,0,t;\frac1{2m})=-S_{+}(\frac12-y,0,t;\frac1{2m})<0$.
Indeed, this follows from   $S_{+}( y,0,t;\frac1{2m})=-\frac1{2m}\frac{\mathrm d}{\mathrm dy} K_{\mathbb R/\mathbb Z}(y,0,t;\frac1{2m})$,
 and the monotonicity properties we have discussed.
\snewpage

Now,
\[\int_{t=0}^\tau\int_{x=0}^1K_{\mathbb R/\mathbb Z}(x,0,t;\tfrac1{2m})S_{+}(y,0,t;\tfrac1{2m})\,\mathrm dx\,\mathrm dt=
\int_{t=0}^\tau S_{+}(y,0,t;\tfrac1{2m})\, \mathrm dt=\]
\[=\sum_{n\in2\mathbb Z}\frac12 \left(\sgn\left( 2\,y+n \right)
- \erf  \left( \frac {\left( 2 y+n \right) \sqrt {m}}{\sqrt {8\tau}}\right)\right)
\]
\[=\sum_{n\in2\mathbb N}\int_{x\in[n+2y,n+2-2y]} \,{\frac {\sqrt {m} }{\sqrt {8\pi\tau } }{{\rm e}^{-\,{\frac {m\,{x}^{2}}{8\tau}}}}}\,\mathrm dx.
\]
This converges to $\frac12-y $ as $\tau\rightarrow+\infty$.
Thus,
\[\int_{t=0}^{+\infty}\int_{x=0}^1K_{\mathbb R/\mathbb Z}(x,0,t;\tfrac1{2m})S_{+}( y,0,t;\tfrac1{2m})\,\mathrm dx\,\mathrm dt=\frac{1 }2-y.\]
Similarly,
\[\int_{t=0}^\tau\int_{x=0}^1-K_{\mathbb R/\mathbb Z}(x,\tfrac12,t;\tfrac1{2m})S_{-}( y,0,t;\tfrac1{2m})\,\mathrm dx\,\mathrm dt
=\int_{t=0}^\tau -S_{-}( y,0,t;\tfrac1{2m})\, \mathrm dt=\]
\[=-\sum_{n\in2\mathbb Z+1}
\frac12\, \left(
\sgn\left( 2\,y+n \right)
- {\erf} \left( {
\frac {\left( 2\,y+n \right) \sqrt {m}}{\sqrt {8\tau}}}
\right)\right)
\]
\[=\sum_{n\in2\mathbb N}\int_{x\in[n+1-2y,n+1+2y]}
 \,{\frac {\sqrt {m} }{\sqrt {8\pi\tau } }{{\rm e}^{-\,{\frac {m\,{x}^{2}}{8\tau}}}}}\,\mathrm dx.
\]
This converges to $y$ as $\tau\rightarrow+\infty$.
Thus,
\[\int_{t=0}^{+\infty}\int_{x=0}^1-K_{\mathbb R/\mathbb Z}(x,\tfrac12,t;\tfrac1{2m})S_{-}( y,0,t;\tfrac1{2m})\,\mathrm dx\,\mathrm dt=y.\]
Therefore, we find
\begin{equation*}
\int_{t\in(0,+\infty)}\int_{x\in(0,1)}  \|A_2^{\gen}(t,x)\|\,\mathrm dx\,\mathrm dt
\leq\frac12\|[Y_1,Y_2]\|,
\end{equation*}
and
\begin{equation*}
\int_{t\in(0,+\infty)}\int_{x\in(0,1)}  A_2^{\gen}(t,x)\,\mathrm dx\,\mathrm dt
= \left(\frac12-2y\right)[Y_1,Y_2].
\end{equation*}

More generally, one can check that the mass primarily generated from
 $ T\cdot \alpha=  T\cdot Y_1\cdot \boldsymbol\delta_{y_1}+ T\cdot Y_2\cdot \boldsymbol\delta_{y_2}$ is according to
\begin{equation*}
\int_{t\in(0,+\infty)}\int_{x\in(0,1)}  A_2^{\gen}(t,x)\,\mathrm dx\,\mathrm dt
= \left(\frac12+y_1-y_2\right)[Y_1,Y_2];
\end{equation*}
with overall variation bounded from $\frac12\|[Y_1,Y_2]\|$.
Furthermore, the overall boundary flux generated  from $\boldsymbol \delta_y\cdot Y$ ($y\in(0,1)$) is  $(\frac12 -y)\cdot Y$
is distributed by a probability measure.

\snewpage

Applied in the formal case $T\cdot\alpha$, we find that the formal solution exist,
 the corresponding version of \eqref{eq:normest1} holds, and $T=1$ can be put the if $\int\|\alpha\|_{\mathfrak g}<1$,
 yielding the corresponding version of Theorem \ref{thm:ex}.
But \eqref{eq:magnustie} is not obtained, as the conjugation is not yet dealt.

The formal conjugating measure to be exponentiated in a time-ordered manner, has density function
\[F^{(T)}(t)=\sum_{n=1}^\infty T^n\cdot F_n(t)\]
with the property
\[\int \|F_1(t)\|\,\mathrm dt\leq   \int_{x\in(0,1)}  \left|\frac12-x\right|\cdot \|\alpha(x)\| \]
for $n=1$, and
\[\int \|F_n(t)\|\,\mathrm dt\leq  \int_{t\in(0,+\infty)}\int_{x\in(0,1)}  \left|\frac12-x\right|\cdot \|A_n^{\gen}(t,x)\|\,\mathrm dx\,\mathrm dt\]
for $n\geq 2$.
Then
\[F_\infty(T) =\exp_{\mathrm R}(t\in (0,+\infty)\mapsto F^{(T)}(t))\]
and
\[H(T)= T\cdot\int_{x\in(0,1)}  \alpha(x)\,\mathrm dx+\sum_{n=2}^\infty T^n\cdot\int_{t\in(0,+\infty)}\int_{x\in(0,1)}  A_n^{\gen}(t,x)\,\mathrm dx\,\mathrm dt\]
has the property
\[\exp(T\cdot \alpha) =F_\infty(T)\cdot (\exp H(T))\cdot (F_\infty(T) )^{-1}
=\exp\left( F_\infty(T)\cdot H(T) \cdot(F_\infty(T) )^{-1} \right).\]
Again, applied in the case when $\alpha$  is a formal noncommutative mass, we see that
\[F_\infty(T) \cdot H(T) \cdot(F_\infty(T) )^{-1}=
\exp_{\mathrm R}(t\in(0,+\infty)\mapsto \ad F^{(T)}(t)) H(T) \]
(strictly speaking only the RHS can be used in a purely Banach--Lie algebraic setting)
must yield the formal Magnus expansion (containing the formal parameter $T$).
But then, again, $T=1$ can be substituted if $\int\|\alpha\|_{\mathfrak g}<1$.
\medskip

\noindent\textit{Sketch of improvement.}
The mass primarily created from  $ T\cdot \alpha=  T\cdot Y_1\cdot \boldsymbol\delta_{y_1}+ T\cdot Y_2\cdot \boldsymbol\delta_{y_2}$
 can be estimated in variation by $J(\dist_{\mathbb R/\mathbb Z} (y_1,y_2))\frac12\|[Y_1,Y_2]\|$ such that
 $\dist_{\mathbb R/\mathbb Z}(y_1,y_2)=\min(\lfloor y_1-y_2\rfloor, \lfloor y_2-y_1\rfloor )$, and
 $J(s)<1$ if $s>0$.
 Taken $y_1,y_2$ randomly, the average estimate is by $ \Delta\cdot\frac12\|[Y_1,Y_2]\|$, where $\Delta<1$.
(Some upper estimate $\Delta$ can be cooked up just by taking $m=1$; and considering a fixed interval in $t$ and an
interval for $\dist_{\mathbb R/\mathbb Z}(y_1,y_2) $.)
If we apply the heat approach for the uniformly distributed formal noncommutative mass in the periodic setting,
 then the mass remains uniformly distributed, and mass creation always occurs in uniformly distributed distances.
Thus \eqref{eq:normest0} can be replaced by
\begin{equation*}
f(T) \stackrel{\forall T^n}{\leq} T\cdot\left(\int\|\alpha\|_{\mathfrak g} \right)+  \Delta\cdot \frac12\cdot\frac12 f(T)^2.
\end{equation*}
Then, ultimately, the bound $\int\|\alpha\|_{\mathfrak g}<1$ can be replaced by $\int\|\alpha\|_{\mathfrak g}<\frac1\Delta$.
\qedremark

\medskip
\noindent\textit{Sketch of improvement.}
In the case of variable diffusion rates, the mass primarily created from  $ T\cdot \alpha=  T\cdot Y_1\cdot \boldsymbol\delta_{y_1}+ T\cdot Y_2\cdot \boldsymbol\delta_{y_2}$ can be estimated in variation by $C\|[Y_1,Y_2]\|$,
where $C$ is an universal constant. ($C=1$ will probably do, but we do claim this.)
However, if the ratio of the diffusion masses tends to infinity, then the total variation of generated mass tends to
$\left(\frac12-\dist_{\mathbb R/\mathbb Z}(y_1,y_2)\right)\|[Y_1,Y_2]\|+\frac14 \|[Y_1,Y_2]\|$
(the first summand coming approximately from the quick diffusion first, the second coming from the slow diffusion later).
Our improvement will apply in the variable-rate case in the setting in the case of formal noncommutative mass
with uniform distribution when $\beta\rightarrow-\infty$.
Then, the variation of mass generated form grades $n'=1$ and $n\geq 2$ can be estimated in average by $(\frac14+\frac14)\|[Y_1,Y_2]\|$.
The variation of the mass generated from   $2\leq n'\leq n$ can be estimated   by $\frac14\|[Y_1,Y_2]\|$
 (a point mass gets introduced into a uniform distribution, essentially).
This means that the limit
\[f(t)=\lim_{\beta\rightarrow -\infty} f^{[\beta]}(T)\]
will be majorized by the solution $g(T)$ of $g(0)=0$ and
\begin{multline*}
g(T)=T\cdot\left(\int\|\alpha\|_{\mathfrak g}\right)+\Delta\,\frac12\cdot T^2\cdot\frac12\left(\int\|\alpha\|_{\mathfrak g}\right)^2 +\\
+ \left(\frac14+\frac14\right)
\cdot T\cdot\left(\int\|\alpha\|_{\mathfrak g}\right)\left( g(T)-T\cdot\left(\int\|\alpha\|_{\mathfrak g}\right)\right)
+\frac14\cdot\frac12\left( g(T)-T\cdot\left(\int\|\alpha\|_{\mathfrak g}\right)\right)^2.
\end{multline*}
Taking $1$ instead of $\Delta$ leads to the criterion in \eqref{eq:magextra}, but we could have done better.

Note  that $g(T)$ does not translate into an upper estimate for the Magnus expansion,
 for that we should involve the terms coming from the extra conjugation.
\qedremark
\begin{remark}
\plabel{rem:higher}
In this section all mass generation estimates were based on primarily created mass from two infinitesimal masses,
 and they can predictably be improved  by using some direct estimates for masses created in higher order.
These are, however, not simple.
\qedremark
\end{remark}

\snewpage

\section{An example in the periodic case}\plabel{sec:perex}
\textbf{The case of precessing measures for $2\times2$ real matrices.}
Here we consider the simplest possible setting where the measures (density functions) are not uniform (constant).
For reasons of tradition, we will use the domain $D=\mathbb R/\pi\mathbb Z$, but this makes no essential difference.
Here we will examine two different questions.
The first one is whether the Maurer-Cartan-heat argument applied to the density $\alpha$
 gives an $F_\infty$ and $H$ such that $\exp_{\mathrm R}(\alpha)=\exp ( F_\infty\cdot H\cdot (F_\infty)^{-1})$.
This is the question of the existence of a ``heat sum'' $F_\infty\cdot H\cdot (F_\infty)^{-1}$.
The second question is whether the heat sum is the same as the sum of the Magnus expansion.
(In what follows the convergence of the Magnus expansion will always be understood as absolute convergence of the Magnus series.)

Now, if $M$ is a $2\times2$ real matrix, then the associated precessing measure is given by
\[ \underbrace{\begin{bmatrix} \cos x & -\sin x\\ \sin x&\cos x\end{bmatrix}
M
\begin{bmatrix} \cos x & -\sin x\\ \sin x&\cos x\end{bmatrix}^{-1} }_{\alpha(x)}\,\mathrm dx|_{[0,\pi]}.
\]
One can see that $\alpha$ can be interpreted as a smooth density function periodic with $\pi$.
It is easy see that the associated time ordered exponential is
\[\exp_{\mathrm R}(\alpha(x)\,\mathrm dx|_{[0,\pi]})=
\exp\left( \pi M+ \pi \begin{bmatrix}&-1\\1&\end{bmatrix} \right)\exp\left( - \pi \begin{bmatrix}&-1\\1&\end{bmatrix} \right) .\]
(Here $\pi$ can be replaced by any positive number.)

We want to solve \eqref{eq:hmc1}--\eqref{eq:hmc2} in the periodic domain.
It is sufficient to consider the case
\[M=\begin{bmatrix} a_0 & c_0\\ b_0& -a_0\end{bmatrix}.\]
Indeed, by conjugation with a rotation matrix $\begin{bmatrix} \cos \xi & -\sin \xi\\ \sin \xi&\cos \xi\end{bmatrix}$,
 which is an inert process here, this shape can be achieved.
Then, a classical solution may be sought in form
\[ A(t,x)=\begin{bmatrix} \cos x & -\sin x\\ \sin x&\cos x\end{bmatrix}
 \begin{bmatrix} a_0 & c(t)\\ b(t)& -a_0\end{bmatrix}
\begin{bmatrix} \cos x & -\sin x\\ \sin x&\cos x\end{bmatrix}^{-1}.
\]
Then $(b(t),c(t))$ must be a solution of the initial value problem
\[b'(t)= k\cdot(-2 )(b(t)+1)(b(t)+c(t)),\]
\[c'(t)= k\cdot  2(c(t)-1)(b(t)+c(t))\]
with
\[b(0)=b_0,\]
\[c(0)=c_0.\]
\snewpage
Regarding this differential equation, it is autonomous, and $(b(t)-1)(c(t)+1) $ turns out to be constant, thus a conserved quantity.
In particular, all trajectories trace parts of hyperbolas.
Another natural symmetry is implemented by $(b(t),c(t))\rightsquigarrow(-c(t),-b(t))$.
The phase diagram of the differential equation is given below.
\[\includegraphics[width=3.5in]{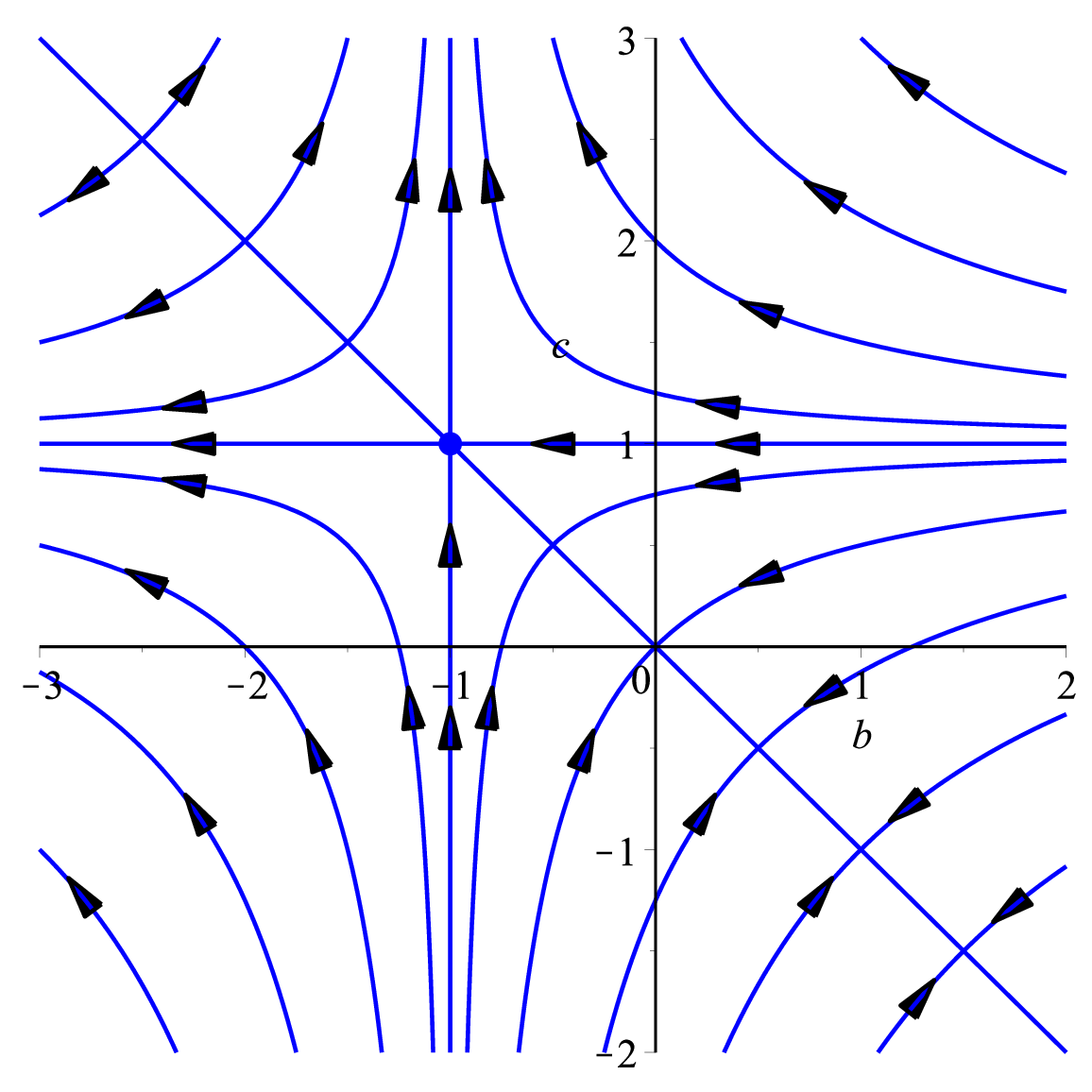}\]
The following table lists all full trajectories up to time translation; $p$ is a positive parameter, $r$ is a real parameter.
\[\begin{array}{c|c|c}
(b(t),c(t))&\text{domain}&(b(t)+1)(c(t)-1)\\
\hline
(-1+\sqrt p \tanh(2k\sqrt pt),1-\sqrt p \coth(2k\sqrt pt) )&t\in(0,+\infty)&-p\\
(-1+\sqrt p \coth(2k\sqrt pt),1-\sqrt p \tanh(2k\sqrt pt) )&t\in(0,+\infty)&-p\\
(-1+\sqrt p \tanh(2k\sqrt pt),1-\sqrt p \coth(2k\sqrt pt) )&t\in(-\infty,0)&-p\\
(-1+\sqrt p \coth(2k\sqrt pt),1-\sqrt p \tanh(2k\sqrt pt) )&t\in(-\infty,0)&-p\\
(-1-\sqrt p \tan(2k\sqrt pt),1-\sqrt p \cot(2k\sqrt pt) )&t\in(0,\pi/2)& p\\
(-1+\sqrt p \cot(2k\sqrt pt),1+\sqrt p \tan(2k\sqrt pt) )&t\in(0,\pi/2)& p\\
(-1+r,1-r )&t\in\mathbb R& -r^2\\
(-1+\frac1{2kt},1 )& t\in(0,+\infty)& 0\\
(-1,1-\frac1{2kt} )&  t\in(0,+\infty)& 0\\
(-1+\frac1{2kt},1 )& t\in(-\infty,0)& 0\\
(-1,1-\frac1{2kt} )&  t\in(-\infty,0)& 0\\
\end{array}\]
\snewpage

We see the following cases:

\textbf{(o)}, where $c_0=-b_0$. In this case the solutions are constant; the heat argument and Magnus expansion works out for obvious reasons.

\textbf{(i)}, where $c_0\neq-b_0$, and $\max (-1 -b_0,c_0-1)=0$.
This is the non-stationary semistable case.
Here the solution exists for infinite time.
Yet it results infinite mass production, the heat sum argument does not work out
due to the problem in the conjugating integral.

Indeed, let us consider the initial condition $b_0= 1, c_0=1$.
Then the solution is given by
\[(b(t),c(t))=\left(-1+\frac1{2(kt+\frac14)},1\right)=\left( \frac{1-4kt}{1+4kt},1\right).\]
Then the flux at the boundary is
\[B(t,0)=k\left.\frac{\partial A(t,x)}{\partial x}\right|_{x=0}=
\begin{bmatrix}
- \frac{2k}{1+4kt}&\\&\frac{2k}{1+4kt}
\end{bmatrix}.\]
This is not Lebesgue integrable for $t\in(0,+\infty)$, in particular, the mass generation is infinite.
This still would not be fatal, but integrated up,
\[F_\tau= \exp_{\mathrm R}( t\in(0,\tau)\mapsto B(t,0))=\begin{bmatrix}
\frac1{\sqrt{1+4k\tau}}&\\&\sqrt{1+4k\tau}
\end{bmatrix}\]
does not converge as $\tau\rightarrow+\infty$.

In fact, we know surely that neither the heat sum nor the Magnus sum can  exist in the case $b_0= 1, c_0=1$.
Indeed, in these cases, the time-ordered exponential
\[
\exp\left(\pi \begin{bmatrix} a_0 & 1\\ 1& a_0\end{bmatrix} +\pi \begin{bmatrix}  & -1\\ 1& \end{bmatrix}\right)
\exp\left(-\pi \begin{bmatrix}  & -1\\ 1& \end{bmatrix}\right)
=-\mathrm e^{a_0}\begin{bmatrix} 1 & \\ 2\pi& 1\end{bmatrix}
\]
would be an exponential of a $2\times2$ real matrix, which is not (cf. \cite{L2}).

Both mechanisms are valid more generally when $b_0\in(-1,+\infty), c_0=1$ or
$b_0=-1, c_0\in(-\infty,1)$.
In these cases,
\begin{equation}
-\mathrm e^{a_0}\begin{bmatrix} 1 & \\ (b_0+1)\pi& 1\end{bmatrix}
\qquad\text{and}\qquad
-\mathrm e^{a_0}\begin{bmatrix} 1&\pi(c_0-1)\\&1\end{bmatrix},
\plabel{eq:exe}
\end{equation}
respectively, will not be an exponentials.

\textbf{(ii)}, where $c_0\neq-b_0$, and $\max (-1-b_0,c_0-1)>0$.
This is the non-stationary  unstable case.
The solution blows up in finite time, and there is no chance for obtaining a heat sum.
Again, we can show that, in this case, the Magnus expansion will not converge.

Indeed, let us first examine the case $(b_0+1)(c_0-1)>0$ (the ``unstable hyperbolic domain'').
Then the time-ordered exponential will give
\[-\mathrm e^{a_0}\begin{bmatrix}
\cosh\sqrt{\pi^2(b_0+1)(c_0-1) } &
\dfrac{ c_0-1}{b_0+1}\sinh\sqrt{\pi^2(b_0+1)(c_0-1) }
\\
\dfrac{b_0+1}{ c_0-1}\sinh\sqrt{\pi^2(b_0+1)(c_0-1) }
& \cosh\sqrt{\pi^2(b_0+1)(c_0-1) } &\end{bmatrix}.\]
This, again, cannot be an exponential, thus, in particular, the Magnus expansion cannot converge.

If $b_0=-1$ or $c_0=1$ but $\max (-1-b_0,c_0-1)>0$ (``the unstable parabolic domain''), then the time-ordered exponentials, as in
 \eqref{eq:exe}, will not be exponentials, again.

Finally, if $b_0<-1$ and $c_0>1$ hold together yet $b_0\neq -c_0$ (``the unstable elliptic domain''), then multiplying the precessing measure
 by an appropriate $\tau\in(0,1)$, we get to the ``unstable hyperbolic domain'' (or just to the ``the semistable parabolic   domain''
 or just to the ``the unstable parabolic  domain'') , where
 the Magnus expansion does not converge.
This precludes the convergence of the original precessing measure, too.

\textbf{(iii)}, where $c_0\neq-b_0$, and $\max (-1-b_0,c_0-1)<0$.
This is the non-stationary  stable case.
Here the solution exists for infinite time.
It converges in exponential rate, which result finite mass production, the argument for heat sum works out.
Does this mean that the Magnus expansion converges?
Not necessarily.
(We are not in the setting of the formal noncommutative masses!)
We claim that that the Magnus expansion will not converge if
\[\min (1-b_0,c_0+1)\leq 0.\]
Indeed, the (absolute) convergence of the Magnus expansion of the density $\alpha(x)$ is equivalent to
 the (absolute) convergence of the Magnus expansion of the density $-\alpha(x)$, and by the previous argument we know
 the negative answer for that.
Thus, here the heat sum provides ``false positives'' for the Magnus expansion.

What finally remains is the case
\begin{equation}
\max(|b_0|,|c_0|)<1.
\plabel{eq:sels}
\end{equation}
Then the Magnus expansion does converge:
Indeed, the result of Moan, Niesen \cite{MN} and Casas \cite{Ca} (about convergence in terms of cumulative operator norm) here guarantees
 the convergence for the Magnus expansion  for $\pi\max(|b_0|,|c_0|)<\pi$;
 and the Magnus sum is the same as the logarithm of the time-ordered exponential.
Using our earlier argument,  we have a direct construction to a solution to the noncommutative heat
 equation under the condition
\begin{equation}
\pi\max(|b_0|,|c_0|)<\frac12.
\plabel{eq:eclas}
\end{equation}
Due to symmetry reasons, the direct construction must also be through precessing measures,
 and then it can be shown to be the same as classical solution we discuss in this section.
By our earlier results, under the condition \eqref{eq:eclas}, the heat sum must be equal to the Magnus sum.
Then, arguments involving real analytical continuation prove that the Magnus sum
 must be equal to the heat sum even under the more general condition \eqref{eq:sels}.
\snewpage

\end{document}